\documentclass{sn-jnl}

\jyear{2021}%

\usepackage{epstopdf}
\usepackage[english]{babel}
\usepackage[T1]{fontenc}
\usepackage[utf8]{inputenc}
\usepackage{calc}

\usepackage{multicol} 
\usepackage{color}

\usepackage{tikz}

\usepackage{calrsfs}
\DeclareMathAlphabet{\pazocal}{OMS}{zplm}{m}{n}

\newsavebox\CBox

\theoremstyle{thmstyleone}%
\newtheorem{theorem}{Theorem}
\newtheorem{proposition}[theorem]{Proposition}%

\theoremstyle{thmstyletwo}%

\theoremstyle{thmstylethree}%
\newcommand\hcancel[2][0.5pt]{%
  \ifmmode\sbox\CBox{$#2$}\else\sbox\CBox{#2}\fi%
  \makebox[0pt][l]{\usebox\CBox}%
  \rule[0.5\ht\CBox-#1/2]{\wd\CBox}{#1}}
  
\newcommand{\R}{\mathbb{R}}

\newcommand{\N}{\mathbb{N}}
\newcommand{\Z}{\mathbb{Z}}

\renewcommand{\P}{\mathbb{P}}
\usepackage{bbm}
\newcommand{\1}{\mathbbm{1}}

\newcommand{\rd}{\mathrm{d}}

\renewcommand{\a}{\alpha}
\renewcommand{\b}{\beta}
\newcommand{\g}{\gamma}

\renewcommand{\d}{\delta}

\newcommand{\e}{\varepsilon}

\renewcommand{\l}{\lambda}

\newcommand{\s}{\sigma}

\renewcommand{\phi}{\varphi}

\newcommand{\noi}{\noindent}

\usepackage{mathtools}

\usepackage{longtable}

\usepackage[subnum]{cases}

\usepackage{natbib}

\usepackage{float}

\makeatletter
\newcommand{\vast}{\bBigg@{3}}
\newcommand{\Vast}{\bBigg@{5}}
\makeatother

\usepackage[title]{appendix}

\usepackage{xr}
\begin{document}

\title[Extremal characteristics of conditional models]{Extremal characteristics of conditional models}


\author*[1]{\fnm{Stan} \sur{Tendijck}}\email{s.tendijck@lancaster.ac.uk}

\author[1]{\fnm{Jonathan} \sur{Tawn}}\email{j.tawn@lancaster.ac.uk}

\author[1,2]{\fnm{Philip} \sur{Jonathan}}\email{p.jonathan@lancaster.ac.uk}

\affil*[1]{\orgdiv{ Department of Mathematics and Statistics}, \orgname{Lancaster University}, \orgaddress{\city{Lancaster}, \postcode{LA1 4YW}, \country{United Kingdom}}}

\affil[2]{\orgname{Shell Research Limited}, \orgaddress{\city{London}, \postcode{SE1 7NA}, \country{United Kingdom}}}


\abstract{Conditionally specified models are often used to describe complex multivariate data. Such models assume implicit structures on the extremes.  So far, no methodology exists for calculating extremal characteristics of conditional models since the copula and marginals are not expressed in closed forms. We consider bivariate conditional models that specify the distribution of $X$ and the distribution of $Y$ conditional on $X$. We provide tools to quantify implicit assumptions on the extremes of this class of models. In particular, these tools allow us to approximate the distribution of the tail of $Y$ and the coefficient of asymptotic independence $\eta$ in closed forms. We apply these methods to a widely used conditional model for wave height and wave period. Moreover, we introduce a new condition on the parameter space for the conditional extremes model of \citet{heffernan2004conditional}, and prove that the conditional extremes model does not capture $\eta$, when $\eta<1$.}


\maketitle

\section{Introduction}
Extreme value theory is a topic of growing interest because of its many important applications in for example risk management \citep{embrechts1999} or ocean engineering \citep{castillo2005extreme}. For instance, in the design or assessment of offshore facilities it is crucial to understand the distribution of extreme sea states. Such extreme sea states are quantified in terms of extreme wave heights, wave periods possibly associated with resonant frequencies, and extreme wind speeds. In risk management, it is important to identify which stocks are likely to suffer extreme losses simultaneously, and to which extent this might happen. In general, we need to use well-estabilished extreme value methods to model such events. Traditionally, such multivariate extreme value methods are composed of marginal models and a dependence copula, each having parametric forms for the tails.

In other areas of statistics, however, it is common to use conditional models for high-dimensional data. Intuitively, this is the most sensible approach. We observe $X$ that partially explains $Y$. So, we define a model for $X$ and a model for $Y$ conditional on $X$. There exist many examples in the literature of models within this conditional framework with applications in extremes, e.g., the conditional extreme value model \citep{heffernan2004conditional, fougeres2012}, the Weibull-log normal distribution \citep[henceforth the Haver-Winterstein distribution]{haver2008environmental}, and hierarchical models \citep{eastoe2019nonstationarity}. Although conditional models are easy to interpret, it can be rather difficult to study the extremes of both $Y$ and $(X,Y)$ within this class. Recently, \citet{engelke2020graphical} developed graphical models for extremes. However, we do not know of any literature that links existing conditional models directly to extremal dependence measures.


There are two extremal dependence measures that are key in identifying and measuring the degree of asymptotic dependence or asymptotic independence \citep{coles1999dependence}. Identifying the correct asymptotic dependence class is important since extrapolation of models from different classes is different. To define asymptotic dependence, we first define $\chi\in[0,1]$, with
\begin{equation}\label{chi_def}
\chi :=\lim_{p\uparrow1}\chi(p) := \lim_{p\uparrow1} \mathbb{P}\left\{Y>F_Y^{-1}(p)\mid X>F_X^{-1}(p)\right\},
\end{equation}
where $F_X$ and $F_Y$ denote the marginal distribution functions of $X$ and $Y$. We say that these random variables are asymptotically dependent if $\chi>0$, i.e., when the joint probability that both random variables are large is of the same magnitude as when one is large. If the coefficient of asymptotic dependence $\chi=0$, we say that the variables are asymptotically independent. In this case, $\chi$ does not give us information on the level of asymptotic independence. So, we additionally define the coefficient of asymptotic independence $\eta\in(0,1]$ \citep{ledford1996statistics} to satisfy for $u\to\infty$
\begin{equation}\label{etaOGdefined}
\mathbb{P}\left\{X>F_X^{-1}\left[F_E(u)\right],\ Y>F_Y^{-1}\left[F_E(u)\right]\right\} \sim \pazocal{L}\left(e^u\right) e^{-u/\eta},
\end{equation}
where $F_E(u) = 1 - \exp(-u)$ is the distribution function of a standard exponential, and where $\mathcal{L}$ is a slowly varying function. Here, we write $f(x)\sim g(x)$ as $x\to\infty$ when $f(x)/g(x)\to1$ as $x\to\infty$. We rewrite definition~(\ref{etaOGdefined}) as
\begin{equation}\label{eta_def}
\eta := \lim_{p\uparrow 1} \eta(p) := \lim_{p\uparrow 1}\frac{\log(1-p)}{\log\left[(1-p)\chi(p)\right]}.
\end{equation}
If the variables are asymptotically dependent, then $\eta=1$; if the variables are asymptotically independent, then $\eta\in(0,1)$ or $\eta=1$ and $\pazocal{L}(u)\to 0$ as $u\to\infty$. The coefficient of asymptotic independence $\eta$ describes the rate of decay to zero of the joint exceedance probability $\mathbb{P}\{X>F_X^{-1}(p),\ Y>F_{Y}^{-1}(p)\}$ as $p$ tends to $1$, see \citet{ledford1996statistics}.

It is relatively straightforward to calculate the two extremal dependence measures for distributions when the joint distribution function is specified parametrically, e.g., a bivariate extreme value distribution \citep{ledford1996statistics}, or when the joint density function is specified parametrically \citep{nolde2021}, e.g., a multivariate normal distribution. In this paper, we consider models specified within the conditional framework. For these cases, it is not straightforward to calculate $\eta$ analytically, and numerical estimation can be difficult since convergence of $\eta(p)$ to $\eta$ can be exceptionally slow. 
We set up methodology to calculate $\eta$ in closed form within this framework and demonstrate the techniques on two widely used examples specified below. We support these limiting results using numerical integration.


First, we consider the model described in \citet{haver2008environmental}, used to explain the dependence between extreme significant wave height and their associated wave periods. 
Secondly, we investigate the model of \citet{heffernan2004conditional}. This is a conditional model which describes the distribution of $Y\mid X$ for large $X$, where both $X$ and $Y$ are on standard margins. As the Heffernan-Tawn model focusses on the averages and deviations of $Y\mid X$ for large $X$, and not necessarily on the tails of $Y\mid X$ for large $X$, it cannot be expected to model $\eta$ correctly. Indeed, we will show that $\eta$ of $(X,Y)$ can be different to $\eta$ from the associated exact Heffernan-Tawn model. More theoretical examples, like $Y\mid X := X^{\beta} Z$ and $Y\mid X := \vert Z\vert^{\vert X\vert}$ where $Z$ is some random variable independent of $X$, can be found in the Ph.D. thesis of \citet{tendijck2023}. 


The layout of the article is as follows. In Section \ref{methodology}, we demonstrate novel techniques for calculating the coefficient of asymptotic independence $\eta$ and illustrate the techniques with some examples. In Sections~\ref{HW_WLN_dist} and~\ref{HTModelCaseStudy}, we apply these techniques to the Haver-Winterstein model and the Heffernan-Tawn model, respectively. Proofs are found in the Appendix and Supplementary Material.

\section{Methodology}\label{methodology}
\subsection{Motivation}

We aim to investigate the extremal properties of the bivariate distribution of $(X,Y)$, for which the distribution of $X$ and the distribution of $Y\mid X$ are specified. In particular, we aim to investigate the tail of the distribution of $Y$ and joint extremes of $X$ and $Y$ via the coefficient of asymptotic independence $\eta$. Deriving such extremal quantities in closed form within this class is not trivial. In this section, we provide a set of tools, derived from the Laplace approximation, to calculate such properties for any conditional model.

First, we consider the tail of the distribution of $Y$. Because the distributions of $X$ and $Y\mid X$ are specified, it is natural to write
\[
1 - F_Y(y) := \P(Y>y) = \int_{-\infty}^{\infty} \P(Y>y\mid X=x)f_X(x)\,\rd x,
\]
where $f_X$ is the density of $X$. In general, this integral is analytically intractable. In Section~\ref{laplace_approx_extension}, we present the tools with which we can derive the asymptotic properties of this integral as $y$ tends to the upper end point of the distribution of $Y$.

To derive the coefficient of asymptotic independence, we additionally need the inverse distribution $F_Y^{-1}(p)$ for values of $p$ close to $1$, and
\[
\P(X>F_X^{-1}(p),\ Y>F_Y^{-1}(p)) = \int_{F_X^{-1}(p)}^{\infty} \P(Y>F_Y^{-1}(p)\mid X = x) f_X(x)\,\rd x.
\]
This integral is also intractable in general; the tools from Section~\ref{laplace_approx_extension} can again be applied to derive the asymptotic decay to $0$ as $p$ tends to $1$.


\subsection{Extension to the Laplace approximation}\label{laplace_approx_extension}
Here we present our theory to calculate asymptotic rates of decay of integrals, that can be used to compute extremal properties, such as $\eta$, of conditional models. We first recall the Laplace approximation, a technique commonly used in Bayesian inference for approximating intractable integrals. This asymptotic approximation forms the basis of our main result. We then state that result, and illustrate key differences with the Laplace approximation by comparing examples.

\begin{proposition}[Laplace approximation]\label{laplace_approximation}
Let $a<b$. Suppose $g:[a,b]\to\R$ is twice continuously differentiable and assume there exists a unique $x^*\in(a,b)$ such that $g(x^*) = \max_{x\in[a,b]}g(x)$ and $g''(x^*)<0$. Then
\[
\int_a^b e^{n g(x) - n g(x^*)}\,\rd x \cdot \sqrt{n(-g''(x^*))} \sim \sqrt{2\pi}
\]
as $n\to\infty$.
\end{proposition}

\noi The main disadvantage of the Laplace approximation is that it can only be used to approximate integrals where the integrands are of the form $f(x)^n$, where $f(x)=e^{g(x)}$ is a positive function. However, we are interested in calculating integrals with integrand $f_n(x) = e^{g_n(x)}$, for some sequence of functions $\{g_n\}_{n\in\N}$.
Now we extend the Laplace approximation under the assumptions that:
(i) the analogue $x_n^*$ of $x^*$ is allowed to depend on $n$; (ii) $x^*_n$ can be equal to either $a$ or $b$; (iii) $g_n''(x^*_n)$ does not need to be negative. 

\begin{proposition}\label{laplace_approximation_variant4}
Let $I\subseteq\R$ be connected with non-zero Lebesgue mass, $k_0\geq 1$ an integer, and $g_n\in C^{k_0}(I)$ a sequence of real-valued (at least) $k_0$-times continuously differentiable functions defined on $I$. For $1\leq i \leq k_0$, we define $g_n^{(i)}$ as the $i$th derivative of $g_n$. We assume that for all $n\in\N$, there exists a unique $x_n^*\in I$ such that $g_n(x_n^*)>g_n(x)$ for all $x\in I\setminus\{x_n^*\}$. Moreover, we assume that $k_0$ is the smallest integer such that $g_n^{(k_0)}(x_n^*) < 0$ and $\lim_{n\to\infty} g_n^{(i)}(x_n^*)[-g_n^{(k_0)}(x_n^*)]^{-i/k_0}=0$ for all $1\leq i<k_0$. Additionally, assume that there exists a $\d$ such that for all $\vert x \vert<\d$
\[
\lim_{n\to\infty} \frac{g_n^{(k_0)}\left\{x_n^* + x\left[-g_n^{(k_0)}(x_n^*)\right]^{-\frac{1}{k_0}}\right\}}{g_n^{(k_0)}(x_n^*)} < \frac{3}{2}.
\]
Then, for $n>N$, there exists a constant $C_1>0$ such that 
\[
\int_{I} e^{g_n(x) - g_n(x_n^*)}\,\rd x \cdot \left[-g_n^{(k_0)}(x_n^*)\right]^{\frac{1}{k_0}} \geq C_1.
\]
\end{proposition}

\noi The proof of Proposition~\ref{laplace_approximation_variant4} can be found in Appendix~\ref{proofs}. One disadvantage of our extension is that it only gives an asymptotic lower bound. In many practical applications, an upper bound can be found directly using inequalities like that in equation~(\ref{upperboundExample}).

\subsection{Examples}
We demonstrate the use of Proposition~\ref{laplace_approximation_variant4} in three cases. Firstly, let $g_n(x)=-n x^p$ for $n\in\N$, $p\in\Z_{\geq1}$ and $I=[0,\infty)$. It is then valid to apply Proposition~\ref{laplace_approximation_variant4} with $x^*_n=0$ and $k_0=p$. Applying the proposition yields a constant $C_1>0$ such that as $n\to\infty$
\[
n^{\frac{1}{p}}\int_{0}^{\infty} e^{-n x^p} \,\rd x \geq C_1.
\]
This lower bound is tight as we now verify for $p\geq 2$, since for $p=1$ the statement holds trivially. For $p\geq 2$, we use the variable transformation $y =  n x^{p}$ to give as $n\to\infty$
\[
n^{\frac{1}{p}}\int_{0}^{\infty} e^{-n x^p} \,\rd x = \frac{1}{p} \int_{0}^{\infty} y^{-\frac{1}{p}-1} e^{-y}  \,\rd y = \Gamma\left(\frac{1}{p}+1\right).
\]
After recognizing that the integral over $[0,\infty)$ is equal to half of the integral over $\R$, we see that Proposition~\ref{laplace_approximation} is also applicable, but only in the special case $p=2$. In this case, Proposition~\ref{laplace_approximation} additionally gives as $n\to\infty$
\[
\int_{0}^{\infty} e^{-n x^2} \,\rd x = \frac{1}{2}\int_{-\infty}^{\infty} e^{-n x^2} \,\rd x \sim \frac{\sqrt{\pi}}{2\sqrt{n}}.
\]
Secondly, let $g_n(x)=-x - n x^2$ and $I=[0,\infty)$. Now Proposition~\ref{laplace_approximation} is not applicable since no function $g(x)$ exists for which $g_n(x)=n g(x)$ holds. Note that Proposition~\ref{laplace_approximation_variant4} is also not applicable with $k_0=1$, since $x_n^*$ has to be equal to $0$ and for $x\neq 0$
\[
\lim_{n\to\infty} \frac{g_n'\left(0 + x\cdot n\right)}{g_n'(0)} = \lim_{n\to\infty} 1+2n^2 x =\infty,
\]
contradicting one of the assumptions. Proposition~\ref{laplace_approximation_variant4} is applicable with $k_0=2$, yielding a constant $C_2>0$ such that as $n\to\infty$
\[
\sqrt{n}\int_{-\infty}^{\infty} e^{-x- n x^2} \,\rd x \geq C_2.
\]
Similar to our first example, this lower bound is tight since we can also directly calculate as $n\to\infty$
\[
\sqrt{n}\int_{-\infty}^{\infty} e^{-x- n x^2} \,\rd x = \sqrt{n}\int_{-\infty}^{\infty} e^{-n\left(x + \frac{1}{2n}\right)^2 + \frac{1}{4n}} \,\rd x \sim \sqrt{\pi}. 
\]
Finally, let $\a_n>0$, $\b_n>0$ for $n\in\N$ and assume $\liminf \alpha_n > 0$. Define $g_n(x) = \alpha_n \log x - \beta_n x$. Using an argument similar to that in the second example, we see that Proposition~\ref{laplace_approximation} is not applicable. However Proposition~\ref{laplace_approximation_variant4} is applicable with $k_0=2$, yielding a constant $C_3>0$ such that as $n\to\infty$
\[
\alpha_n^{-\alpha_n - \frac{1}{2}} \beta_n^{\alpha_n+1} e^{\alpha_n} \int_{0}^{\infty} x^{\alpha_n} e^{- \beta_n x} \,\rd x \geq C_3.
\]
This bound is also tight, which can be seen from recognizing the density of a gamma distribution in the expression above, and applying limit results for the gamma function.

\section{Haver-Winterstein model}\label{HW_WLN_dist}
\citet{haver2008environmental} introduce the Haver-Winterstein (HW) distribution for significant wave height $H_S$ and wave period $T_p$ in the North Sea. Their model is set up in the conditional framework: they specify a class of distributions for $H_S$ and a class of distributions for $T_p\mid H_S$. Variations of this approach have been widely applied in ocean engineering with over 150 citations, 25 of which correspond to 2021, see for example \citet{drago2013}. However we are not aware of any literature quantifying $\chi$ and $\eta$ in closed form for the HW distribution; we now show how to calculate these.



The HW distribution is formulated as
\begin{equation}\label{HW_HS}
f_{X}(x) = \begin{cases} \frac{1}{\sqrt{2\pi} \alpha x} \exp\left\{-\frac{(\log x - \theta)^2}{2\alpha^2}\right\},\  & \text{for}\ 0<x\leq u, \\ \frac{k}{\l^{k}} x^{k-1} \exp\left\{-\left(\frac{x}{\l}\right)^{k}\right\},\ & \text{for}\ x>u. \end{cases}
\end{equation}
where $u,\alpha,k,\l>0$ and $\theta\in\R$. In particular, the parameters are constrained such that $f_X$ is continuous at $u$ and integrates to $1$. Secondly, they take $Y\mid X$ to be conditionally log-normal 
\begin{equation}\label{HW_TP}
f_{Y\mid X}(y\mid x) = \frac{1}{\sqrt{2\pi} \sigma(x) y} \exp\left\{-\frac{(\log y - \mu(x))^2}{2\sigma(x)^2}\right\},\ \ \ \ \text{for}\ x,y>0,
\end{equation}
where $\mu(x):=\mu_0+\mu_1 x^{\mu_2}$ and $\sigma(x):=\left[\sigma_0 + \sigma_1 \exp(-\sigma_2 x)\right]^{1/2}$ with $\mu_0\in\R,\ \mu_1,\mu_2,\sigma_0,\sigma_1,\sigma_2 > 0$.



Model parameter estimates \citep{haver2008environmental} from data observed in the northern North Sea are given in the Supplementary Material. For ease of presentation, we make two assumptions about the parameter space of the HW distribution that are consistent with parameter estimates $(\hat{\mu}_2,\hat{k})=(0.225,1.55)$ from \citet{haver2008environmental}. Specifically, we make the following restrictions: $0<\mu_2<0.5$ and $2\mu_2<k$. These assumptions reduce the number of cases to be considered significantly whilst including realistic domains for the parameters as considered by practioners.

We now show how to use Proposition~\ref{laplace_approximation_variant4} to calculate the extremal dependence measures $\chi$ and $\eta$ for the bivariate random vector $(X,Y)$ distributed according to the HW 
distribution in the restricted parameter space. Calculation of $\eta$ is split into two steps. In the first step, we calculate the distribution function $F_Y$ of $Y$ and in the second we evaluate the rate of decay of joint probabilities $\P\{X > F_X^{-1}[F_E(u)],Y>F_Y^{-1}[F_E(u)]\}$ as $u$ tends to infinity.

\begin{figure}
\centering
\includegraphics[width=0.8\textwidth]{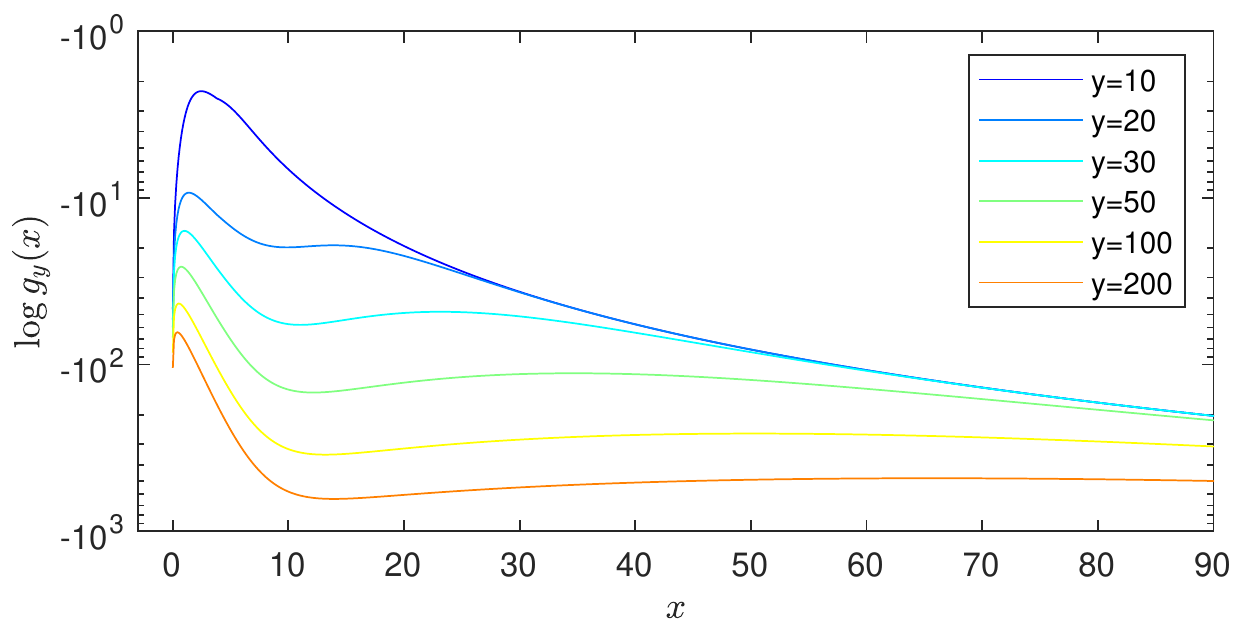}
\caption{The function $\log g_y$ from equation~(\ref{gy}) for $y=10,\ 20,\ 30,\ 40,\ 50,\ 100$ with parameters as reported in~\citet{haver2008environmental}, see Supplementary Material.}
\label{logDensity}
\end{figure}

We have
\begin{equation}\label{WLN_integral1}
\P(Y > y) = \int_0^{\infty} \P(Y> y\mid X=x)f_X(x) \,\rd x = \int_0^{\infty} \overline{\Phi}\left(\frac{\log y - \mu(x)}{\sigma(x)}\right)f_X(x)\,\rd x,
\end{equation}
where $\overline{\Phi}$ is the survival function of a standard Gaussian. This integral is analytically intractable but we can calculate its limiting leading order behaviour in closed form. Proposition~\ref{laplace_approximation_variant4} gives a lower bound and an upper bound of the same order as the lower bound is then found directly. For ease of notation, we denote the integrand by
\begin{equation}\label{gy}
g_y(x) :=  \overline{\Phi}\left(\frac{\log y - \mu(x)}{\sigma(x)}\right)f_X(x)
\end{equation}
for $x>0$. In Figure~\ref{logDensity}, we plot $g_y$ for various values of $y$. From the figure, we note that $g_y$ has two local maxima for suffiiciently large $y$. These are $x_y^*$, which converges to zero, and $x_y^{**}$, which diverges to infinity. This observation implies that we cannot apply Proposition~\ref{laplace_approximation_variant4} directly in this case. We therefore proceed as follows: (i) calculate $x_y^*$ and $x_y^{**}$; (ii) partition the interval of integration into intervals $I_1$ and $I_2$, where $x_y^*\in I_1$ and $x_y^{**}\in I_2$, such that the conditions of Proposition~\ref{laplace_approximation_variant4} hold for both intervals, and then apply the proposition on each interval; (iii) combine the two lower bounds found to get a lower bound for integral~(\ref{WLN_integral1}); (iv) derive a limiting upper bound for integral~(\ref{WLN_integral1}) of the same order as the lower bound.

In the Supplementary Material, we derive that as $y\to\infty$
\[x_y^* \sim \left(\frac{\s_1 \s_2\cdot \log y}{2\mu_1\mu_2(\s_0+\s_1)}\right)^{-\frac{1}{1-\mu_2}}\ \ \ \text{and}\ \ \ x_y^{**} \sim \left(\frac{ \l^k\mu_1\mu_2\cdot \log y}{k \s_0}\right)^{\frac{1}{k-\mu_2}}.\]
From Figure~\ref{logDensity}, we recognize that $g_y(x_y^*) > g_y(x_y^{**})$ as $y\to\infty$. We show that this holds analytically in the Supplementary Material when $2\mu_2<k$. We now apply Proposition~\ref{laplace_approximation_variant4} and find that $k_0=2$ is appropriate. The proposition then gives a lower bound for integral~(\ref{WLN_integral1}) around $x_y^*$ as $y\to\infty$ of
\[
\P(Y>y)  \geq  \exp\left\{-\frac{\log^2 y}{2(\s_0+\s_1)} + O(\log y)\right\}.
\]
Finally, since $g_y(x_y^*)>g_y(x_y^{**})$, it is straightforward to show as $y\to\infty$ that
\[
\P(Y>y) \leq  \exp\left\{-\frac{\log^2 y}{2(\s_0+\s_1)} + O(\log y)\right\}
\]
using the inequality
\begin{equation}\label{upperboundExample}
\P(Y> y\mid X=x)f_X(x) \leq g_y(x_y^*)\1\{x\in[0,x_y^{**}]\} + f_X(x) \1\{x>x_y^{**}\}.
\end{equation}
We now can calculate $\eta$ and show that $\chi=0$. To that end, we first need to calculate the inverse probability integral transform, transforming $Y$ to standard exponential margins; i.e., we need $F_Y^{-1}[F_E(u)]$. Next, we need to evaluate the asymptotic behaviour of $\P\{Y>F_Y^{-1}[F_E(u)],X>F_X^{-1}[F_E(u)]\}$ as $u\to\infty$. To evaluate $F_Y^{-1}\circ F_E$, we first calculate for $y\to\infty$
\[\begin{aligned}
F_E^{-1}(F_Y(y)) = -\log(1-F_Y(y))= \frac{\log^2 y}{2(\s_0+\s_1)} + O(\log y).
\end{aligned}\]
We invert this expression by solving $F_E^{-1}(F_Y(y))=u$ for $\log y$. This yields $\log y = \sqrt{2\sigma_0+\sigma_1)u} + O(1)$ as $y\to\infty$. We can now write down an asymptotic expression for $\chi(u)$ as $u\to\infty$
\[\begin{aligned}
\chi(u)&:=\P\left\{F_E^{-1}\left[F_Y(Y)\right] > u,\ F_E^{-1}\left[F_X(X)\right] > u\right\} \\
&= \P\left\{\log Y >\sqrt{2}(\sigma_0+\sigma_1) \sqrt{u} + O(1),\ (X/\l)^k > u\right\}  \\
& = \int_{\l u^{1/k}}^{\infty} \overline{\Phi}\left(\frac{\sqrt{2}(\s_0+\s_1)\sqrt{u} + O(1) - \mu(x)}{\sigma(x)}\mid X=x\right)\cdot \frac{k x^{k-1}}{\l^k} \exp\left\{-\left(\frac{x}{\l}\right)^k\right\}\,\rd x.
\end{aligned}\]
In the Supplementary Material, we show that Proposition~\ref{laplace_approximation_variant4} is applicable for this integral with $k_0=1$ and $x_u^*=\l u^{1/k}$. Moreover, we derive directly an upper bound of the same order, obtaining
\[
\chi(u) =\exp\left\{- \left(2  +\frac{\s_1}{\s_0}\right) u + O\left(u^{1/2 + \mu_2/k}\right)\right\}
\]
as $u\to\infty$. Hence, $\chi=0$ and
\[\eta = \left(2  +\frac{\s_1}{\s_0}\right)^{-1}.\]
In particular, for the parameter estimates from \citet{haver2008environmental}, the value of $\eta\in(0,1/2)$ implies that the distribution exhibits negative asymptotic independence \citep{ledford1996statistics}.


\section{Heffernan-Tawn model}\label{HTModelCaseStudy}
In multivariate extreme value theory, the conditional extreme value model of \citet{heffernan2004conditional}, henceforth denoted the HT model, is widely studied and applied to extrapolate multivariate data. The HT model has been cited over 600 times, and is applied e.g. in oceanography \citep{ross2020}, finance \citep{hilal2011}, and spatio-temporal extremes \citep{simpson2021}. The HT model is a limit model and its form is motivated by derived limiting forms from numerous theoretical examples. \citet{keef2013estimation} assume that for $(X,Y)$ on standard Laplace margins there exist parameters $\a\in[-1,1]$, $\b<1$ and a non-degenerate distribution function $H$ such that for $x>0$, $z\in\R$
\begin{equation}\label{htlimiteq}
\lim_{u\to\infty} \P\left(\frac{Y-\a X}{X^{\b}} \leq z,\ X-u>x\mid X>u\right) = \exp(-x) H(z).
\end{equation}
In the limit of $u\to\infty$, this formulation implies that $(Y-\a X)X^{-\b}$ and $(X-u)$ are independent conditional on $X>u$, and are distributed as $H$ and a standard exponential, respectively. As is common practice in extreme value theory, the HT model assumes that the corresponding limiting family in~(\ref{htlimiteq}) holds exactly at a finite level. Thus the HT model is specified for $x>u$, where $u$ is a sufficiently high threshold such that the limit representation in~(\ref{htlimiteq}) is considered a good approximation. Let $(X,Y)$ be a random vector such that $X$ and $Y$ both have standard Laplace margins. Moreover, let $\a,\b\in[0,1)$ and assume that for $x>u>0$
\begin{equation}\label{HT_formulation}
\P(Y> y\mid X=x) = \overline{H}\left(\frac{y - \a x}{x^{\b}}\right)
\end{equation}
holds for all $y\in\R$ where $\overline{H} = 1- H$ is some non-degenerate survival function. In this case, we say that $(X,Y)$ are distributed according to an exact version of the HT model. We consider two cases for $H$, corresponding to finite and infinite upper end points. If $H$ has a finite upper end point $z^H$, calculations for $\eta$ are trivial. Indeed, when $X=x$, $Y$ cannot be larger than $\a x + x^{\b} z^H$. In particular, as $u\to\infty$, $Y>u$ is equivalent to $X>u/\alpha+O(u^{\b})$. Hence as $u\to\infty$
\[\begin{aligned}
\P(X>u,Y>u) &\sim \P\left\{X>u,X>u/\alpha+O(u^{\b})\right\} \\
&\sim \P\left\{X>u/\alpha+O(u^{\b})\right\} \\
&=\exp\left\{-u/\alpha+O(u^{\b})\right\}.
\end{aligned}\]
Therefore, $\eta=\a$ when $\a>0$ and otherwise does not exist. 

Now assume that $H$ has an infinite upper end point. To make calculations tractable, we parameterise $\overline{H}$ as
\begin{equation}\label{H_parameterisation}
\overline{H}\left(z\right) = \exp\left\{- \gamma z^{\delta} + o\left(z^{\delta}\right)\right\}\1\{z>0\} + \1\{z\leq 0\}
\end{equation}
for $\gamma>0$, $\delta\geq 1$. For simplicity, we do not consider potential negative arguments for $\overline{H}$ since the precise form of its lower tail is not relevant to the current work. Parameterisation~(\ref{H_parameterisation}) covers most non-trivial cases for the upper tail including Gaussian, Weibull and exponential tails; see examples in \citet{heffernan2004conditional}. Moreover if the tail of $\overline{H}$ is heavier than that of the exponential, $Y$ cannot possibly possibly follow a standard Laplace distribution. This links to the restricton $\d\geq 1$. For illustration, we set $ o(z^{\delta})=0$ in equation~(\ref{H_parameterisation}). The resulting Weibull survival function is a suitable choice for $\overline{H}$, since it has an extreme value tail index of $0$, but a varying tail thickness controlled by $\delta$.

\begin{proposition}\label{lower_bound_delta}
If $(X,Y)$ follows distribution~(\ref{HT_formulation}) with $H$ as in~(\ref{H_parameterisation}) with $o(z^{\d})=0$, then $\d\geq(1-\beta)^{-1}$.
\end{proposition}

\noi The proof of Proposition~\ref{lower_bound_delta} is found in Appendix~\ref{proofs}. Following similar arguments to those used in the proof of Proposition~\ref{lower_bound_delta}, we calculate $\chi$ and $\eta$ for any combination of the parameters $(\alpha,\beta,\delta,\gamma)$ in their specified parameter space. We collect results in  Table~\ref{table_HT_prm}. In the Supplementary Material, we only give details of the $\eta$ calculations when $\a,\b\in(0,1)$, $\g>0$ and $\d=(1-\b)^{-1}$. For the other five cases in Table~\ref{table_HT_prm}, we state results without proof. In particular, the argument underpinning the $\eta$ calculation when $\d>(1-\b)^{-1}$ is similar to the argument used when $\overline{H}$ has a finite upper end point. In this case, $\eta=\a$ when $\a>0$ and when $\a=0$, $\eta$ is not defined.

In Table~\ref{table_HT_prm}, it is convenient to refer to $c=\max\{1,c_0\}\in[1,1/\alpha)$ where $c_0\in(0,1/\alpha)$ satisfies
\begin{equation}\label{c0}
\gamma(1-\alpha c_0)^{\delta-1}\left(\d-1 +\alpha c_0\right) = c_0^{\delta}.
\end{equation}

\begin{table}
\centering
\begin{tabular}{c|c|c|c|c}
$\alpha$ & $\beta$ & $\gamma$ & $\delta$ & $\eta$ \\ \hline
$(0,1)$ & $[0,1)$ & $(0,\infty)$ & $\left((1-\b)^{-1},\infty\right)$ & $\a$ \\\hline
$(0,1)$ & $(0,1)$ & $(0,\infty)$ & $(1-\b)^{-1}$ & $\left(\frac{\g(1-\alpha c)^{\delta}}{c^{\delta-1}} + c\right)^{-1}$ \\\hline
$(0,1)$ & $0$ & $(1/\alpha,\infty)$ & $1$ & $\alpha$ \\\hline
$(0,1)$ & $0$ & $(0,1/\alpha]$ & $1$ & $1/(\g+1-\g\a)$ \\\hline
$0$ & $(0,1)$ & $(0,\infty)$ & $\left((1-\b)^{-1},\infty\right)$ & Not defined \\\hline
$0$ & $(0,1)$ &  $(0, (1-\beta)/\beta]$  & $(1-\b)^{-1}$ & $1/(\g+1)$ \\\hline
$0$ & $(0,1)$ &$[(1-\beta)/\beta,\infty)$& $(1-\b)^{-1}$ & $\g^{-1/\d} (\d-1)^{1-1/\d}/\d$ \\
\end{tabular}
\caption{Values of $\eta$ for model~(\ref{HT_formulation}) with $\overline{H}$ as in~(\ref{H_parameterisation}) for different ranges of parameter combinations, where $c=\max\{1,c_0\}\in[1,1/\alpha)$ for $c_0$ given in equation~(\ref{c0}).}
\label{table_HT_prm}
\end{table}

\begin{figure}
\centering
\includegraphics[width=\textwidth]{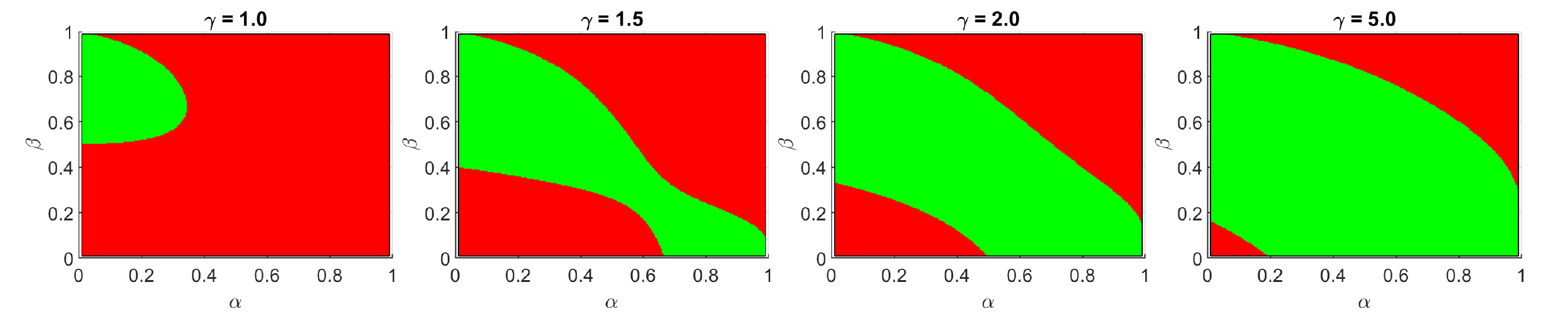}
\caption{Visualisation of $c_0$ from equation~(\ref{c0}) for $\gamma=1,\ 1.5,\ 2,\ 5$ and $\d=(1-\b)^{-1}$. The region corresponding to $c_0\in(0,1)$ is shown in red; the region corresponding to $c_0\in(1,1/\alpha)$ is shown in green.}
\label{cEqualToOne}
\end{figure}

\begin{figure}
\centering
\includegraphics[width=\textwidth]{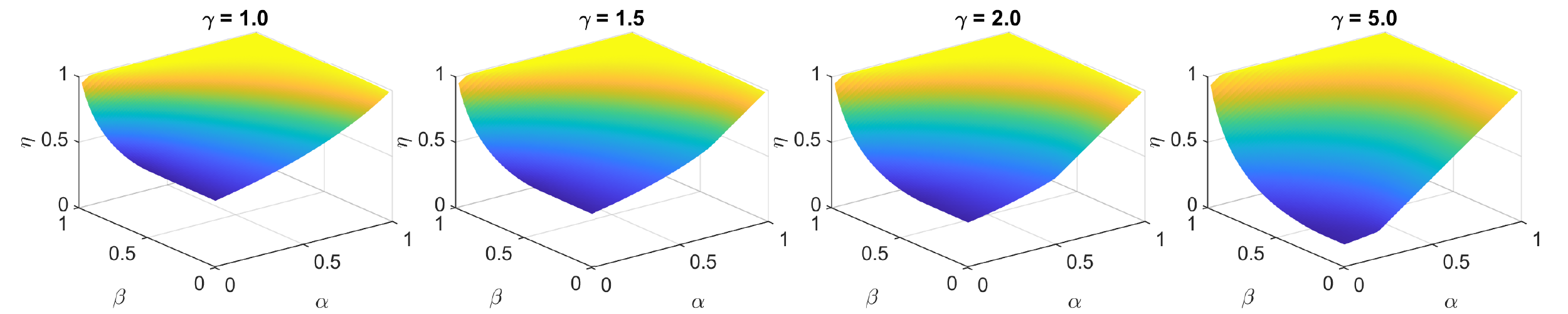}
\caption{The value of $\eta$ as a function of $\alpha$, $\beta$ and $\gamma$ with $\d=(1-\b)^{-1}$ from the HT model~(\ref{HT_formulation}) and~(\ref{H_parameterisation}).}
\label{htModelPrmComb}
\end{figure}

\noi To give some intuition on the value of $c$, in Figure~\ref{cEqualToOne} we sketch the region of the parameter space corresponding to $c=1$ (in red) for different values of $\gamma$. Finally in Figure~\ref{htModelPrmComb} we visualise $\eta$ for a set of different parameter combinations with $\d=(1-\b)^{-1}$. 

We note the following interesting findings. The parameter $\eta$ is non-decreasing with increasing $\alpha$ and with increasing $\beta$. Parameter combinations $(\a,\b,\g,\d)$ exist for which $\a,\b>0$ but $\eta<0.5$. Hence, there are cases for which $Y$ increases with $X$ but the extremes of $(X,Y)$ are negatively associated as measured by $\eta$ \citep{ledford1996statistics}.

Finally we note that the Heffernan-Tawn model is not $\eta$ invariant, i.e., when the HT model occurs in the limit of the distribution of $(X,Y)$, then $\eta$ for $(X,Y)$ is not necessarily the same as $\eta$ for the associated exact HT model. To illustrate this, let $(X,Y)$ follow an inverted bivariate extreme value distribution with a logistic dependence structure \citep{ledford1996statistics} on Laplace margins with parameter $\xi\in(0,1]$, such that
\begin{equation}\label{inv_log}
\P(X>x,\ Y>y) = \exp\left\{-\left[t_x^{1/\xi} + t_{y}^{1/\xi}\right]^{\xi}\right\},
\end{equation}
where $t_x := \log 2 - \log[2-\exp(x)]$ for $x<0$ and $t_x:=\log 2 + x$ for $x>0$, with $t_y$ similarly defined. It is straightforward to derive that in the limit, the Heffernan-Tawn model~(\ref{HT_formulation}) is applicable to $(X,Y)$ with $\overline{H}$ as in equation~(\ref{H_parameterisation}) and $o(z^{\d})=0$. Specifically, 
\[
\lim_{x\to\infty}\P\left(Y X^{\xi - 1}> z \mid X=x\right)  = \exp\left(- \xi z^{1/\xi}\right).
\]
Now let $(X_{HT},Y_{HT})$ be distributed following our exact version of the HT model associated with $(X,Y)$. That is,  for $X_{HT}<u$, we have $(X_{HT},Y_{HT})=(X,Y)$. For $X_{HT}\geq u$, $X_{HT}-u$ is standard exponentially distributed, and $Y_{HT}\mid X_{HT}$ follows model~(\ref{HT_formulation}) with $\overline{H}$ as in~(\ref{H_parameterisation}) with parameters $(\alpha,\ \beta,\ \gamma,\ \delta)=(0,\ 1-\xi,\ \xi,\ 1/\xi)$ and $o(z^{\d})=0$. In this case $\gamma < (1-\beta)/\beta$, and Table~\ref{table_HT_prm} implies that the coefficient of asymptotic independence $\eta_{HT}$ of $(X_{HT},Y_{HT})$ is equal to $1/(\xi+1)$. In contrast, it is straightforward to derive directly from definition~(\ref{inv_log}) that $\eta$ of $(X,Y)$ is equal to $2^{-\xi}$. Hence $\eta_{HT}\neq \eta$ when $\xi\in(0,1)$. 

Finally we illustrate numerically the differences between $\eta$, $\eta_{HT}$ and their finite level counterparts $\eta(p)$ and $\eta_{HT}(p)$ for $p\in(0,1)$. 
For definiteness, we let $(X,Y)$ follow distribution~(\ref{inv_log}) with $\xi=0.35$. We simulate a sample $\{(x_i,y_i):\ i=1,\dots,n\}$ of size $n=10,000$. First we empirically estimate $\eta(p)$ from equation~(\ref{eta_def}) for $p\in(0,1)$ and calculate pointwise $95\%$ confidence intervals using the binomial distribution. Next we note that $\eta(p)=\eta$ for $p\in(0.5,1)$. Finally we calculate the corresponding $\eta_{HT}(p)$ for $p$ near $1$ using numerical integration. 

\begin{figure}
\centering
\includegraphics[width=0.9\textwidth]{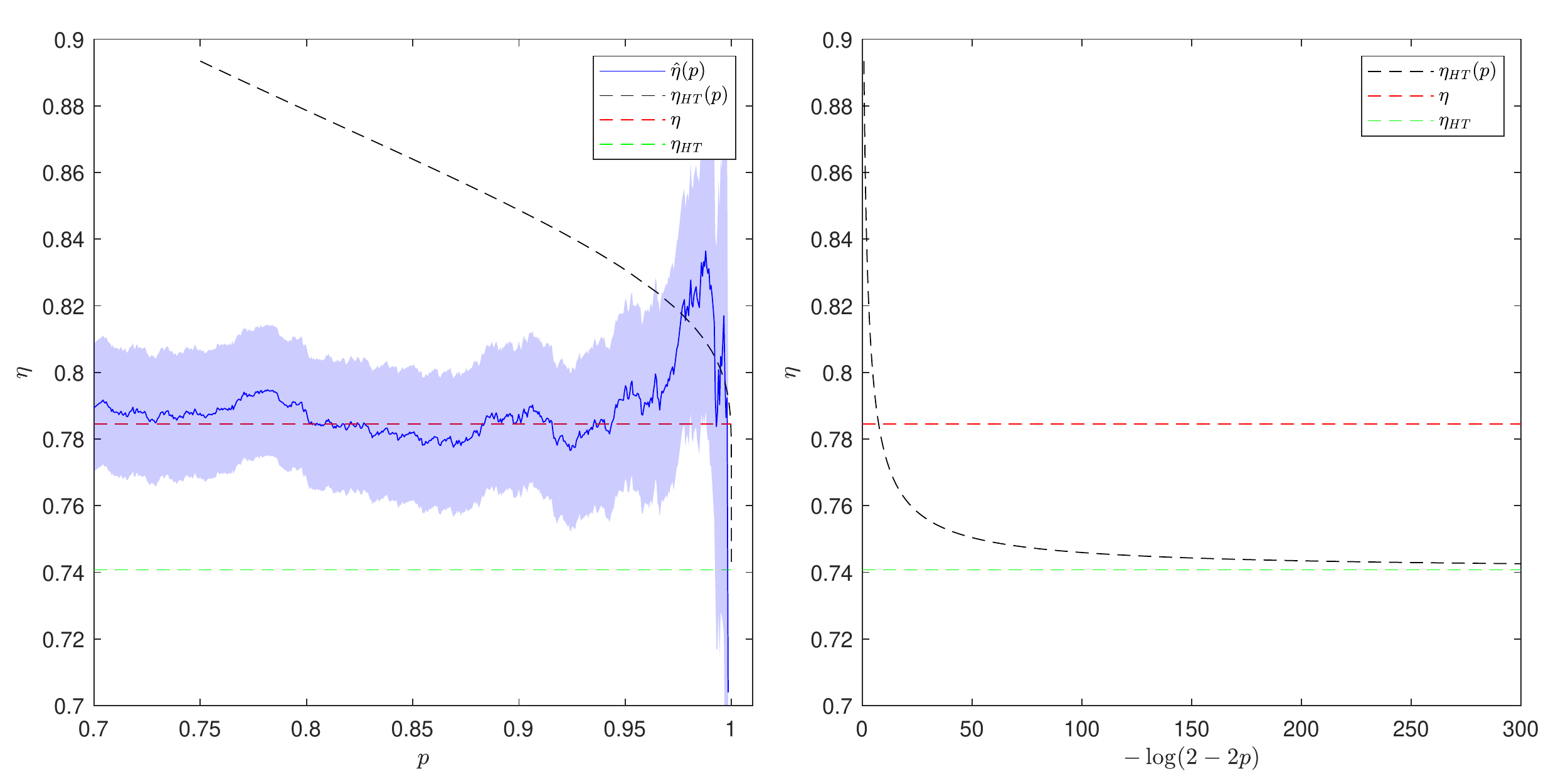}
\caption{Coefficients of asymptotic independence $\eta$ (red dashed) for distribution~(\ref{inv_log}) with $\xi=0.35$, and the corresponding value for the exact limiting HT model $\eta_{HT}$ (green dashed), and its finite level counterpart $\eta_{HT}(p)$ (black dashed). Empirical estimates $\hat{\eta}(p)$ for a sample of size $10,000$ with pointwise confidence intervals are shown in blue. Left and right hand panels are the same except for the scale of the $x$-axis, set on the right to illustrate the behaviour of $\eta_{HT}(p)$ for $p$ near $1$.}
\label{logisticEta}
\end{figure}

Results are shown in Figure~\ref{logisticEta}. Left and right hand plots are the same except for the scale of the $x$-axis, illustrating the behaviour of $\eta_{HT}(p)$ for $p$ near $1$. Reassuringly, the true $\eta$ of the underlying model (red dashed) falls within the $95\%$ confidence interval for its empirical counterpart $\hat{\eta}(p)$ (blue). Further, $\eta_{HT}(p)$ (black dashed) converges to $\eta_{HT}$ (green dashed). We note that $\eta_{HT}(p)$ varies as a function of $p$ and only seems to asymptote for $p>1-\exp(-50)/2 \approx 1 - 9.6\cdot 10^{-23}$. Finally, since $\eta_{HT}<\eta$, we would expect that $\eta_{HT}(p)$ would underestimate $\eta$, but it turns out this is only the case for $p>1-\exp(-7.5)/2\approx 0.9997$.

\bmhead{Supplementary information}
In the Supplementary Material, we give details of the mathematical derivations corresponding to the case studies.

\bmhead{Acknowledgments}

We acknowledge motivating discussions with Ed Mackay of Exeter University and David Randell of Shell. This article is based on work completed while Stan Tendijck was part of the EPSRC funded STOR-i centre for doctoral training (grant no. EP/L015692/1), with part-funding from Shell Research Ltd.

\begin{appendices}

\section{Proofs}\label{proofs}
\noi \emph{Proof of Proposition~\ref{laplace_approximation_variant4}.} 
We prove that as $n\to\infty$, there exists a constant $C_1>0$ such that
\[
\mathcal{I}_n := \int_{I} e^{g_n(x) - g_n(x_n^*)}\,\rd x \cdot \left(-g_n^{(k_0)}(x_n^*)\right)^{\frac{1}{k_0}} \geq C_1.
\]
To bound $\mathcal{I}_n$ from below, we first simplify its expression by applying the variable transformation $y = t_n(x):= (x-x_n^*)\left(-g_n^{(k_0)}(x_n^*)\right)^{1/k_0}$ and defining
\[
h_n(y) := g_n\left(x_n^* + y \left(-g_n^{(k_0)}(x_n^*)\right)^{-\frac{1}{k_0}}\right), \ \ \text{for}\ y \in I_n' := \left\{t_n(x):\ x\in I\right\}.
\]
Then, the integral $\mathcal{I}_n$ becomes
\[\begin{aligned}
\mathcal{I}_n = \int_{I_n'} e^{h_n(y)- h_n(0)}\,\rd y.
\end{aligned}\]
We note that for all $n\in\N$, we have $0\in I_n'$, $h_n\in C^{k_0}(I_n')$, and $h_n(0)>h_n(y)$ for all $y\in I_n'\setminus\{0\}$. Moreover, we have for $y\in I_n'$, $i=1,\dots,k_0$,
\[
h_n^{(i)}(y) =  g^{(i)}_n\left(x_n^* + y\left(-g_n^{(k_0)}(x_n^*)\right)^{-1/k}\right)\cdot \left(-g_n^{(k_0)}(x_n^*)\right)^{-i/k_0}.
\]
Hence, $h_n^{(k_0)}(0)=-1$ and $\lim_{n\to\infty} h_n^{(i)}(0)=0$ for all $1\leq i<k_0$. Using Taylor's theorem, there exists a function $\xi(y)$ taking on a value between $0$ and $y$ such that
\[\begin{aligned}
h_n(y) - h_n(0) & = \sum_{i=1}^{k_0-1} \frac{y^i}{i!} h_n^{(i)}(0)+ \frac{y^{k_0}}{k_0!} h_n^{(k_0)}(\xi(y)).
\end{aligned} \]
Let $\e>0$. Because $\lim_{n\to\infty} h_n^{(i)}(0)=0$ for all $i<k_0$, we can find an $N_0\in\N$ such that for all $n>N_0$, we have $\max_{i=1,\dots,k_0-1}|h_n^{(i)}(0)|<\e$. Moreover, from the assumptions of the proposition, we can find a $\d>0$ and an $N_1\in\N$ such that for all $n>N_1$, $h_n^{(k_0)}(y)>-3/2$ for $y\in(-\d,\d)\cap I_n'$. For $n>\max\{N_0,N_1\}$, 
\[\begin{aligned}
h_n(y) - h_n(0) 
& > -|y| \epsilon - \frac{|y|^2}{2!} \epsilon - \cdots - \frac{|y|^{k_0-1}}{(k_0-1)!} \e - \frac{3|y|^{k_0}}{2k_0!}  > -\e e^{\d}- \frac{3|y|^{k_0}}{2k_0!} 
\end{aligned} \]
for $y \in (-\d,\d)\cap I_n'$. Hence, we derive a lower bound
\[
\mathcal{I}_n \geq e^{ -\e e^{\d}} \int_{I_n' \cap (-\d,\d)} e^{-\frac{3 |y|^{k_0}}{2 k_0!}}\,\rd y =: C_1.
\]
From the connectedness of $I$ and $0\in I_n'$, we conclude that $I_n' \cap (-\d,\d)$ has positive mass under the Lebesgue measure. Hence, $C_1\in(0,\infty)$.
 \hfill $\square$ \medskip

\noi \emph{Proof of Proposition~\ref{lower_bound_delta}.}
Let $(X,Y)$ be a random vector such that $X$ and $Y$ both have standard Laplace margins. Moreover, assume that there exist $-1\leq \a\leq 1$, $0\leq \b<1$ and $u>0$ such that for $x>u$
\[
\P(Y> y\mid X=x) = \overline{H}\left(\frac{y - \a x}{x^{\b}}\right)
\]
holds for all $y\in\R$ with
\[
\overline{H}\left(z\right) = \exp(- \gamma z^{\delta})\1\{z>0\} + \1\{z\leq 0\},
\]
where $\gamma,\delta>0$. We now derive that $\delta \geq (1-\beta)^{-1}$ must hold. Since $Y$ is distributed as a standard Laplace, we have for $y>0$
\[\begin{aligned}
\frac{\exp(-y)}{2}
& = \P(\a X + X^{\b} Z \geq y,\ X\geq u) \P(X\geq u) + \P(Y \geq y,\ X< u) \\
& \geq \P(\a X + X^{\b} Z \geq y,\ X\geq u)   \geq \P(X^{\b} Z \geq y,\ X\geq u)  \\
& = \int_{u}^{\infty} \P\left(Z \geq \frac{y}{x^{\b}}\right)f_X(x)\,\rd x   = \frac{1}{2}\int_{u}^{\infty} \exp\left(- \frac{\gamma y^{\d}}{x^{\b\d}} - x\right)\,\rd x  =: \tilde{\mathcal{I}}_y. \\
\end{aligned}\]
We will show that $2\exp(y)\tilde{\mathcal{I}}_y > 1$ as $y\to\infty$ if $\delta < (1-\beta)^{-1}$, which thus would contradict with the marginal distribution of $Y$. This result holds trivially for $\beta=0$. So, for now, we let $\beta>0$. We will prove this asymptotic inequality by applying Proposition~\ref{laplace_approximation_variant4}, with $k_0=2$, to bound $\tilde{\mathcal{I}}_y$ from below. 

First define $I:=[u,\infty)$ as the integration domain, and
\[
g_y(x):= \exp\left(- \frac{\gamma y^{\d}}{x^{\b\d}} - x\right)\1\{x\in I\},\ \ \ \text{and}\ \ \ 
h_y(x):=\left(- \frac{\gamma y^{\d}}{x^{\b\d}} - x\right)\1\{x\in I\}.
\]
Next we find the mode $x_y^*$ of $g_y(x)$. We assume that $x_y^*$ lies in the interior of $I$ such that $h_y'(x_y^*)=0$, which implies that $\b\d \gamma y^{\d} (x_y^*)^{-\b\d-1} = 1$. So, $x_y^* = \left(\b\d \gamma\right)^{\frac{1}{\b\d+1}} y^{\frac{\d}{\b\d+1}}$,  which lies in the interior of $I$ for sufficiently large $y$. We now compute
\[
g_y(x_y^*) =  \exp\left(- \frac{\gamma y^{\d}}{(x_y^*)^{\b\d}} - x_y^*\right)
=  \exp\left(- A y^{\frac{\d}{\b\d+1}}\right)
\]
with $A:=\gamma\left(\b\d \gamma\right)^{-\frac{\b\d}{\b\d+1}} +  \left(\b\d \gamma\right)^{\frac{1}{\b\d+1}}$.
Secondly,
\[
h_y''(x_y^*) =  - \b\d(\b\d+1) (x_y^*)^{-\b\d-2} \gamma y^{\d}
=  - (\b\d+1)\left(\b\d\g\right)^{-\frac{1}{\b\d+1}} y^{-\frac{\d}{\b\d+1}}.
\]
Using these expression, we can now check that the assumptions from Proposition~\ref{laplace_approximation_variant4} with $k_0=2$ are satisfied. First we note that $h_y'(x_y^*)(-h_y''(x_y^*))^{-1/2} = 0$. Next let $C>0$ and $|x|\leq C$, then
\[\begin{aligned}
\lim_{y\to\infty} & \frac{h_y''\left(x_y^*+\frac{x}{\sqrt{-h_y''(x_y^*)}}\right)}{h_y''(x_y^*)}  \\
& = \lim_{y\to\infty} \frac{ - \b\d(\b\d+1)\left(\left(\b\d \gamma\right)^{\frac{1}{\b\d+1}} y^{\frac{\d}{\b\d+1}}+\frac{x}{\sqrt{ (\b\d+1)\left(\b\d\g\right)^{-\frac{1}{\b\d+1}} y^{-\frac{\d}{\b\d+1}}}}\right)^{-\b\d-2} \gamma y^{\d}}{ - (\b\d+1)\left(\b\d\g\right)^{-\frac{1}{\b\d+1}} y^{-\frac{\d}{\b\d+1}}} \\
& = \lim_{y\to\infty} \frac{\left(y^{\frac{\d}{\b\d+1}}+\frac{x}{\sqrt{ (\b\d+1)\left(\b\d\g\right)^{\frac{1}{\b\d+1}} y^{-\frac{\d}{\b\d+1}}}}\right)^{-\b\d-2} y^{\d}}{ y^{-\frac{\d}{\b\d+1}}} \\
& = \lim_{y\to\infty} \left(1+\frac{x}{\sqrt{ (\b\d+1)\left(\b\d\g\right)^{\frac{1}{\b\d+1}} y^{\frac{\d}{\b\d+1}}}}\right)^{-\b\d-2} \\
& = 1,
\end{aligned}\]
which is sufficient to show that for each $\tilde{x}$, Proposition~\ref{laplace_approximation_variant4} is applicable with $k_0=2$ on interval $I_{\tilde{x}} := \left[x_y^*-\frac{\tilde{x}}{\sqrt{-h_y''(x_y^*)}},x_y^*+\frac{\tilde{x}}{\sqrt{-h_y''(x_y^*)}}\right]$. Hence for each $\tilde{x}$, there exists a constant $C_1(\tilde{x})>0$ such that as $y\to\infty$
\[\begin{aligned}
y^{-\frac{\d/2}{\b\d+1}} \exp\left(A y^{\frac{\d}{\b\d+1}}\right)\cdot \tilde{\mathcal{I}}_y
&\geq y^{-\frac{\d/2}{\b\d+1}} \exp\left(A y^{\frac{\d}{\b\d+1}}\right)\cdot\int_{I_{\tilde{x}}} g_y(x)\,\rd x \\
&= \frac{C_1(\tilde{x})  \left(\b\d\g\right)^{\frac{1}{2(\b\d+1)}} }{\sqrt{\b\d+1}}.
\end{aligned}\]
Using the inequality $2\exp(y) \tilde{\mathcal{I}}_y \leq 1$ as $y\to\infty$, we must have
\begin{equation}\label{ineqProp}
 y^{-\frac{\d/2}{\b\d+1}} \exp\left(A y^{\frac{\d}{\b\d+1}}\right) \cdot \frac{1}{2}\exp(-y) \geq  \frac{C_1(\tilde{x})  \left(\b\d\g\right)^{\frac{1}{2(\b\d+1)}} }{\sqrt{\b\d+1}}
\end{equation}
as $y\to\infty$. Since $0\leq \b<1$, we note that if $\d<(1-\b)^{-1}$ then inequality~(\ref{ineqProp}) does not hold. So, we derive that $\d \geq (1-\b)^{-1}$. \hfill $\square$ \medskip

\end{appendices}


\bibliography{bibpaper3}

\newcommand{\noop}[1]{}
\begin{thebibliography}{}

\bibitem[\protect\astroncite{Castillo et~al.}{2005}]{castillo2005extreme}
Castillo, E., Hadi, A.~S., Balakrishnan, N., and Sarabia, J.-M. (2005).
\newblock Extreme {V}alue and {R}elated {M}odels with {A}pplications in
  {E}ngineering and {S}cience.
\newblock Hoboken, New Jersey: Wiley.

\bibitem[\protect\astroncite{Coles et~al.}{1999}]{coles1999dependence}
Coles, S.~G., Heffernan, J.~E., and Tawn, J.~A. (1999).
\newblock Dependence measures for extreme value analyses.
\newblock {\em Extremes}, 2(4):339--365.

\bibitem[\protect\astroncite{Drago et~al.}{2013}]{drago2013}
Drago, M., Giovanetti, G., and Pizzigalli, C. (2013).
\newblock Assessment of significant wave height--peak period distribution
  considering the wave steepness limit.
\newblock In {\em International Conference on Offshore Mechanics and Arctic
  Engineering}, volume 55393, page V005T06A015. American Society of Mechanical
  Engineers.

\bibitem[\protect\astroncite{Eastoe}{2019}]{eastoe2019nonstationarity}
Eastoe, E.~F. (2019).
\newblock Nonstationarity in peaks-over-threshold river flows: A regional
  random effects model.
\newblock {\em Environmetrics}, 30(5):e2560.

\bibitem[\protect\astroncite{Embrechts et~al.}{1999}]{embrechts1999}
Embrechts, P., Resnick, S.~I., and Samorodnitsky, G. (1999).
\newblock Extreme value theory as a risk management tool.
\newblock {\em North American Actuarial Journal}, 3(2):30--41.

\bibitem[\protect\astroncite{Engelke and Hitz}{2020}]{engelke2020graphical}
Engelke, S. and Hitz, A.~S. (2020).
\newblock Graphical models for extremes (with discussion).
\newblock {\em Journal of the Royal Statistical Society: Series B (Statistical
  Methodology)}, 82(4):871--932.

\bibitem[\protect\astroncite{Fougeres and Soulier}{2012}]{fougeres2012}
Fougeres, A.-L. and Soulier, P. (2012).
\newblock Estimation of conditional laws given an extreme component.
\newblock {\em Extremes}, 15(1):1--34.

\bibitem[\protect\astroncite{Haver and
  Winterstein}{2009}]{haver2008environmental}
Haver, S. and Winterstein, S.~R. (2009).
\newblock Environmental contour lines: A method for estimating long term
  extremes by a short term analysis.
\newblock {\em Transactions of the Society of Naval Architects and Marine
  Engineers}, 116:116--127.

\bibitem[\protect\astroncite{Heffernan and
  Tawn}{2004}]{heffernan2004conditional}
Heffernan, J.~E. and Tawn, J.~A. (2004).
\newblock A conditional approach for multivariate extreme values (with
  discussion).
\newblock {\em Journal of the Royal Statistical Society: Series B
  (Methodology)}, 66(3):497--546.

\bibitem[\protect\astroncite{Hilal et~al.}{2011}]{hilal2011}
Hilal, S., Poon, S.-H., and Tawn, J.~A. (2011).
\newblock Hedging the black swan: Conditional heteroskedasticity and tail
  dependence in {S}\&{P}500 and {V}{I}{X}.
\newblock {\em Journal of Banking \& Finance}, 35(9):2374--2387.

\bibitem[\protect\astroncite{Keef et~al.}{2013}]{keef2013estimation}
Keef, C., Papastathopoulos, I., and Tawn, J.~A. (2013).
\newblock Estimation of the conditional distribution of a multivariate variable
  given that one of its components is large: Additional constraints for the
  {H}effernan and {T}awn model.
\newblock {\em Journal of Multivariate Analysis}, 115:396--404.

\bibitem[\protect\astroncite{Ledford and Tawn}{1996}]{ledford1996statistics}
Ledford, A.~W. and Tawn, J.~A. (1996).
\newblock Statistics for near independence in multivariate extreme values.
\newblock {\em Biometrika}, 83(1):169--187.

\bibitem[\protect\astroncite{Nolde and Wadsworth}{2021}]{nolde2021}
Nolde, N. and Wadsworth, J.~L. (2021).
\newblock Linking representations for multivariate extremes via a limit set.
\newblock {\em Advances in Applied Probability, to appear}.

\bibitem[\protect\astroncite{Ross et~al.}{2020}]{ross2020}
Ross, E., Astrup, O.~C., Bitner-Gregersen, E., Bunn, N., Feld, G., Gouldby, B.,
  Huseby, A., Liu, Y., Randell, D., Vanem, E., et~al. (2020).
\newblock On environmental contours for marine and coastal design.
\newblock {\em Ocean Engineering}, 195:106194.

\bibitem[\protect\astroncite{Simpson and Wadsworth}{2021}]{simpson2021}
Simpson, E.~S. and Wadsworth, J.~L. (2021).
\newblock Conditional modelling of spatio-temporal extremes for {R}ed {S}ea
  surface temperatures.
\newblock {\em Spatial Statistics}, 41:100482.

\bibitem[\protect\astroncite{Tendijck}{2023}]{tendijck2023}
Tendijck, S. H.~A. (2023).
\newblock Multivariate {O}ceanographic {E}xtremes in {T}ime and {S}pace.
\newblock Ph.D. thesis at Lancaster University (United Kingdom), to appear in
  2023.

\end{thebibliography}


 \newcommand{\noop}[1]{}
\begin{thebibliography}{}

\bibitem[\protect\astroncite{Haver and
  Winterstein}{2009}]{haver2008environmental}
Haver, S. and Winterstein, S.~R. (2009).
\newblock Environmental contour lines: A method for estimating long term
  extremes by a short term analysis.
\newblock {\em Transactions of the Society of Naval Architects and Marine
  Engineers}, 116:116--127.

\end{thebibliography}
\bibliographystyle{apa}


\end{document}


\title[Extremal characteristics of conditional models]{Extremal characteristics of conditional models}


\author*[1]{\fnm{Stan} \sur{Tendijck}}\email{s.tendijck@lancaster.ac.uk}

\author[1]{\fnm{Jonathan} \sur{Tawn}}\email{j.tawn@lancaster.ac.uk}

\author[1,2]{\fnm{Philip} \sur{Jonathan}}\email{p.jonathan@lancaster.ac.uk}

\affil*[1]{\orgdiv{ Department of Mathematics and Statistics}, \orgname{Lancaster University}, \orgaddress{\city{Lancaster}, \postcode{LA1 4YW}, \country{United Kingdom}}}

\affil[2]{\orgname{Shell Research Limited}, \orgaddress{\city{London}, \postcode{SE1 7NA}, \country{United Kingdom}}}


\maketitle

\section{Introduction}

We give an overview of the content in the Supplementary Material. In Section \ref{hwmodelprm_supmat}, we give parameter estimates of the Haver-Winterstein (HW) distribution as referred to in Section \ref{HW_WLN_dist}. In Sections~\ref{hwmodeldetails_supmat_1}-\ref{hwmodeldetails_supmat_6}, we give the details of the calculations that support the arguments in Section \ref{HW_WLN_dist}. Finally, in Section~\ref{supHT_section} one can find the mathematical derivations of the results stated in of Section \ref{HTModelCaseStudy}.

\section{Supplementary Material}
\subsection{HW model parameters}\label{hwmodelprm_supmat}

\begin{table}[!htbp]
\centering
\begin{tabular}{|l|c|c|c|c|c|c|} \hline
Parameter & $\alpha$ & $\theta$ & $u$ & $k$ & $\l$ & \\\hline
& 0.573 & 0.893 & 3.803 & 1.550 & 2.908 & \\\hline
Parameter & $\mu_0$& $\mu_1$& $\mu_2$& $\sigma_0$& $\sigma_1$& $\sigma_2$ \\\hline
& 1.134 & 0.892 & 0.225 & 0.005 & 0.120 & 0.455 \\\hline
\end{tabular}
\caption{Parameters of the joint probability density function of significant wave height $H_S$ (m) and wave period $T_p$ (s) for the Northern North Sea \citep{haver2008environmental}.}
\label{HW_par_est}
\end{table}

\subsection{Details on calculations for the HW distribution}\label{hwmodeldetails_supmat_1}
Let $(X,Y)$ follow the HW model, see Section~\ref{HW_WLN_dist}, with $0<\mu_2<0.5$ and $2\mu_2<k$. The goal is to calculate the asymptotic behaviour of joint probabilities $\P(X>F_X^{-1}(p),\ Y>F_Y^{-1}(p))$ when $p$ tends to $1$ where $F_X$ and $F_Y$ denote the distribution functions of the random variables $X$ and $Y$, respectively. First, we evaluate the distribution function of $Y$ at large values such that we can calculate $F_Y^{-1}(p)$. After, we compute joint probabilities, like $\P(X_E > u,Y_E>u)$, where $X_E$ and $Y_E$ denote $X$ and $Y$, respectively, on exponential margins, i.e., $X_E=-\log(1 - F_X(X))$ and $Y_E=-\log(1-F_Y(Y))$.

We write down an analytical expression for the survival function $\overline{F}_Y$ of $Y$
\begin{equation}\label{equation_FY}
\overline{F}_Y(y) := \P(Y>y) = \int_0^{\infty} \overline{\Phi}\left(\frac{\log y - \mu(x)}{\sigma(x)}\right)f_X(x)\,\rd x.
\end{equation}
where $\mu(x)$ and $\sigma(x)$ are defined in the main paper. We remark that we need to evaluate $\overline{F}_Y$ at large $y$. To that end, we denote $p_y(x):=(\log y - \mu(x))/\sigma(x)$, and the integrand
\begin{equation}\label{gy22}
g_y(x) :=  \overline{\Phi}\left(p_y(x)\right)f_X(x).
\end{equation}
As seen in Figure~\ref{logDensity}, the integrand $g_y$ has two local maxima for $y$ large enough. Hence, Proposition~\ref{laplace_approximation_variant4} is not directly applicable. However, we can use the proposition to indirectly prove a lower bound for the intergal~(\ref{equation_FY}). Next, it is straightforward to find an upper bound for the integral with the same rate of decay as the proven lower bound.

We follow the following sets of steps: (a) we prove that there exist (at least) two local maxima $x_y^*$ and $x_y^{**}$, and find expressions for both. If there are more then $x_y^*$ is the one with the smallest argument, and $x_y^{**}$ is the one with the second smallest argument; (b) we show that we can apply part of Proposition~\ref{laplace_approximation_variant4} to the smaller of the two local maxima, which gives a lower bound for the integral; (c) we define an upperbound $\tilde{g}_y$ for the integrand $g_y$, compute the integral of $\tilde{g}_y$, and show that this integral has the same rate of decay as the lower bound; (d) finally, we combine the results to get an asymptotic expression for $\overline{F}_Y(y)$ as $y\to\infty$.

We need to start by working out the expressions for the local maxima. We do this by considering all possible options, which yields five (types of) local extrema $0<x_0<x_1<x_2<x_3<x_4<\infty$ that satisfy the following: (i) as $y\to\infty$, $p_y(x_0) \to \infty$ holds and $x_0\to0$; (ii) as $y\to\infty$, $p_y(x_1) \to \infty$ holds and $x_1\to c\in(0,\infty)$; (iii) as $y\to\infty$, $p_y(x_2) \to \infty$ holds and $x_2\to\infty$; (iv) as $y\to\infty$, $p_y(x_3) \to c\in\R$ holds and $x_3\to\infty$; (v) as $y\to\infty$, $p_y(x_4) \to -\infty\in\R$ holds and $x_4\to\infty$. It is straightforward to show that $x_3$ and $x_4$ cannot exist. However, this argument is unnecessary for the purpose of this section.

Finally, we calculate $\overline{F}_Y(y)$ using Proposition~\ref{laplace_approximation_variant4}. In particular, we will get a lower bound by applying Proposition~\ref{laplace_approximation_variant4} around the local maximum $x_0$ and we derive an upper bound directly.

\subsection{Finding local extrema}
We consider cases (i), (ii) and (iii). These cases have in common that $p_y(x_*) \to \infty$ for $x_*\in\{x_0,x_1,x_2\}$. We will write $x_*$ rather than either $x_0,x_1,x_2$ to not distinguish between arguments that are applicable to all three cases. To find an expression for $x_*$ in closed form, we define $h_y(x):=\log g_y(x)$ and we solve $h_y'(x_*)=0$. First, we calculate $h_y'(x)$,
\[\begin{aligned}
h_y'(x) &= \frac{\rd}{\rd x}\left(\log\overline{\Phi}\left(p_y(x)\right) + \log f_X(x)\right) \\
& = \frac{-\phi\left(p_y(x)\right)}{\overline{\Phi}\left(p_y(x)\right)} \cdot p_y'(x) + \frac{f'_X(x)}{f_X(x)}.
\end{aligned} \]
Since $p_y(x_*)\to\infty$, we can simplify this expression by using Mills' ratio, which says that
\[
\frac{\overline{\Phi}(x)}{\phi(x)} = \frac{1}{x} - \frac{1}{x^3} + O(x^{-5})
\]
as $x\to\infty$, which implies $\phi(x)/\overline{\Phi}(x) = x + x^{-1} + O(x^{-3})$ as $x\to\infty$. Moreover, we can write
\[\begin{aligned}
p_y'(x)&=\frac{\rd}{\rd x}\left(\frac{\log y - \mu(x)}{\sigma(x)}\right) \\
& = -(\log y - \mu(x))\cdot \frac{\sigma'(x)}{\sigma(x)^2} - \frac{\mu'(x)}{\sigma(x)}\\
& =- p_y(x)\cdot\frac{\sigma'(x)}{\sigma(x)} - \frac{\mu'(x)}{\sigma(x)}.
\end{aligned} \]
So,
\[\begin{aligned}
h_y'(x_*) &=-\left(p_y(x_*)+\frac{1}{p_y(x_*)} + O\left(p_y(x_*)^{-3}\right)\right) \cdot \left(- p_y(x_*)\cdot\frac{\sigma'(x_*)}{\sigma(x_*)} - \frac{\mu'(x_*)}{\sigma(x_*)}\right)\\
&\text{\hspace{1cm}}\ + \frac{f'_X(x_*)}{f_X(x_*)} \\
&=p_y(x_*)^2\cdot \frac{\sigma'(x_*)}{\sigma(x_*)} + p_y(x_*) \cdot \frac{\mu'(x_*)}{\sigma(x_*)}+\frac{\sigma'(x_*)}{\sigma(x_*)} + \frac{\mu'(x_*)}{p_y(x_*)\sigma(x_*)} \\
&\text{\hspace{1cm}}\ + O\left(\frac{\sigma'(x_*)}{p_y(x_*)^2 \sigma(x_*)} + \frac{\mu'(x_*)}{p_y(x_*)^3 \sigma(x_*)}\right) + \frac{f'_X(x_*)}{f_X(x_*)} \\
\end{aligned}\]
as $y\to\infty$. We now fill in the parametric forms for $\mu$ and $\sigma$. We can then simplify this expression to

\[ \begin{aligned}
h_y'(x_*) = &\frac{(\log y - \mu_0 - \mu_1 x_*^{\mu_2})^2 \cdot -\frac{1}{2} \s_1\s_2 \exp\left(-\sigma_2 x_*\right)(\s_0+\s_1\exp(-\s_2 x_*))^{-1/2} }{(\s_0+\s_1\exp(-\s_2 x_*))^{3/2}} \\
&+ \frac{(\log y - \mu_0 - \mu_1 x_*^{\mu_2})\mu_1\mu_2 x_*^{\mu_2-1}}{\s_0+\s_1\exp(-\s_2 x_*)} \\
&-\frac{\frac{1}{2}\s_1\s_2 \exp(-\s_2 x_*) (\s_0+\s_1 \exp(-\s_2 x_*))^{-1/2}}{(\s_0+\s_1\exp(-\s_2 x_*))^{1/2}}\\
&+ \frac{\mu_1 \mu_2 x_*^{\mu_2-1}}{\log y - \mu_0 -\mu_1x_*^{\mu_2}} \\
&+ O\left(\frac{x_*^{\mu_2-1}}{(\log y - \mu_0-\mu_1 x_*^{\mu_2})^3} + \frac{\exp(-\s_2 x_*)}{(\log y-\mu_0-\mu_1x_*^{\mu_2})^2}\right) + \frac{f_X'(x_*)}{f_X(x_*)} \\
= &-\frac{(\log y - \mu_0 - \mu_1 x_*^{\mu_2})^2 \cdot  \s_1\s_2 \exp\left(-\sigma_2 x_*\right) }{2(\s_0+\s_1\exp(-\s_2 x_*))^{2}} \\
&+ \frac{(\log y - \mu_0 - \mu_1 x_*^{\mu_2})\mu_1\mu_2 x_*^{\mu_2-1}}{\s_0+\s_1\exp(-\s_2 x_*)} -\frac{\s_1\s_2 \exp(-\s_2 x_*)}{2(\s_0+\s_1\exp(-\s_2 x_*))}\\
&+ \frac{\mu_1 \mu_2 x_*^{\mu_2-1}}{\log y - \mu_0 -\mu_1x_*^{\mu_2}} + O\left(\frac{x_*^{\mu_2-1}}{(\log y - \mu_0-\mu_1 x_*^{\mu_2})^3} + \frac{\exp(-\s_2 x_*)}{(\log y-\mu_0-\mu_1x_*^{\mu_2})^2}\right) \\
&+ \frac{f_X'(x_*)}{f_X(x_*)}. \\
\end{aligned} \]
Since, $h_y'(x_*)=0$ for all $y$, we can let $y\to\infty$, to further simplify
\begin{equation}\label{h_y_prime_x0_linear} \begin{aligned}
0 &= \lim_{y\to\infty} h_y'(x_*) \\
&= \lim_{y\to\infty}\Bigg(-\log^2 y \cdot \frac{\s_1\s_2\exp(-\s_2 x_*)}{2(\s_0+\s_1 \exp(-\s_2 x_*))^2} \\
&\text{\hspace{2cm}}\ + \log y \cdot\left( \frac{(\mu_0 + \mu_1 x_*^{\mu_2}) \s_1\s_2 \exp\left(-\sigma_2 x_*\right)}{(\s_0+\s_1\exp(-\s_2 x_*))^{2}}+  \frac{\mu_1\mu_2 x_*^{\mu_2-1}}{\s_0+\s_1\exp(-\s_2 x_*)} \right) \\
&\text{\hspace{2cm}}\  - \frac{(\mu_0 + \mu_1 x_*^{\mu_2})\mu_1\mu_2 x_*^{\mu_2-1}}{\s_0+\s_1\exp(-\s_2 x_*)} -\frac{\s_1\s_2 \exp(-\s_2 x_*)}{2(\s_0+\s_1\exp(-\s_2 x_*))}  \\
&\text{\hspace{2cm}}\ - \frac{ (\mu_0 + \mu_1 x_*^{\mu_2})^2 \cdot  \s_1\s_2 \exp\left(-\sigma_2 x_*\right) }{2(\s_0+\s_1\exp(-\s_2 x_*))^{2}} + \frac{\mu_1 \mu_2 x_*^{\mu_2-1}}{\log y - \mu_0 -\mu_1x_*^{\mu_2}} + \frac{f_X'(x_*)}{f_X(x_*)}\Bigg).
\end{aligned}\end{equation}
We now split up the analysis into the three cases: (i) $x_*=x_0\to0$; (ii) $x_* = x_1\to c\in(0,\infty)$; (iii) $x_*=x_2\to\infty$.

\bigskip

\noi \textbf{Case (i): $x_*=x_0\to0$}

\bigskip

\noi In this case, there must exist a $y'>0$ such that for all $y>y'$, $x_0(y) < u$. So, let $y>y'$, then
\[
\frac{f_X'(x_0)}{f_X(x_0)} = -\frac{\log x_0 - \theta}{x_0\alpha^2} - \frac{1}{x_0}.
\] 
Filling in $x_*=x_0$ simplifies equation~(\ref{h_y_prime_x0_linear}) to
\begin{equation}\label{case1} \begin{aligned}
\lim_{y\to\infty}&\Bigg(-\log^2 y \cdot \frac{\s_1 \s_2}{2(\s_0+\s_1)^2}+\log y \cdot \left(\frac{\mu_0 \s_1\s_2}{(\s_0+\s_1)^2} + \frac{\mu_1\mu_2 x_0^{\mu_2-1}}{\s_0+\s_1}\right)  - \frac{\mu_0 \mu_1\mu_2 x_0^{\mu_2-1}}{\s_0+\s_1}\\
& -\frac{\s_1\s_2}{2(\s_0+\s_1)} - \frac{\mu_0^2\s_1\s_2}{2(\s_0+\s_1)^2} +\frac{\mu_1 \mu_2 x_0^{\mu_2-1}}{\log y} -\frac{\log x_0 - \theta}{x_0\alpha^2} - \frac{1}{x_0}\Bigg) = 0.
\end{aligned} \end{equation}
Because $0<\mu_2<0.5$, the dominating terms within this limit are of the order $\log^2(y)$ and $\log y \cdot x_0^{\mu_2-1}$. Indeed, $(\log x_0)/x_0$ is dominated by both of these terms since, we must eventually have $x_0^{2\mu_2-2}> (\log x_0)/x_0$. So $x_0$ must satisfy as $y\to\infty$
\[
-\log y \cdot \frac{\s_1 \s_2}{2(\s_0+\s_1)}+\log y \cdot \mu_1\mu_2 \cdot x_0^{\mu_2-1} = O\left(\frac{\log x_0}{x_0 \log y}\right).
\]
Finally, we derive the following asymptotic expression for $x_0$ as $y\to\infty$
\begin{equation}\label{0}
x_0  = \left(\frac{\s_1 \s_2}{2\mu_1\mu_2(\s_0+\s_1)}\right)^{-\frac{1}{1-\mu_2}} \cdot \left(\log y\right)^{-\frac{1}{1-\mu_2}} + O\left(\log^{-2}(y)\right).
\end{equation}
We will later show that $h_y''(x_0) < 0$. So, indeed $x_0$ corresponds to a local maximum.

\bigskip

\noi \textbf{Case (ii): $x_*=x_1\to c\in(0,\infty)$}

\bigskip

\noi In this case, equation~(\ref{h_y_prime_x0_linear}) is equivalent to
\[
 \lim_{y\to\infty}-c_1 \log^2 y + c_2 \log y  - c_3 = 0
\]
where
\[ \begin{aligned}
0 < c_1 &=  \frac{\s_1\s_2\exp(-\s_2 c)}{2(\s_0+\s_1 \exp(-\s_2 c))^2} \\
0 < c_2 &=  \frac{(\mu_0 + \mu_1 c^{\mu_2}) \s_1\s_2 \exp\left(-\sigma_2 c\right)}{(\s_0+\s_1\exp(-\s_2 c))^{2}}+  \frac{\mu_1\mu_2 c^{\mu_2-1}}{\s_0+\s_1\exp(-\s_2 c)} \\
0 < c_3 &= \frac{(\mu_0 + \mu_1 c^{\mu_2})\mu_1\mu_2 c^{\mu_2-1}}{\s_0+\s_1\exp(-\s_2 c)} + \frac{\s_1\s_2 \exp(-\s_2 x_*)}{2(\s_0+\s_1\exp(-\s_2 x_*))}  \\
&\text{\hspace{1cm}}\ + \frac{ (\mu_0 + \mu_1 c^{\mu_2})^2 \cdot  \s_1\s_2 \exp\left(-\sigma_2 c\right) }{2(\s_0+\s_1\exp(-\s_2 c))^{2}} + \frac{\mu_1 \mu_2 c^{\mu_2-1}}{\log y - \mu_0 -\mu_1 c^{\mu_2}} - \frac{f_X'(c)}{f_X(c)}.
\end{aligned}\]
We can now clearly see that equation~(\ref{h_y_prime_x0_linear}) cannot be valid under this assumption. We conclude that $x_1$ does not exist. 

\bigskip

\noi\textbf{Case (iii): $x_*=x_2\to\infty$}

\bigskip

\noi Finally, let $x_*=x_2\to\infty$. In this case, there must exist a $y''>0$ such that for all $y>y''$, $x_0(y) > u$. So, let $y>y''$, then
\[
\frac{f_X'(x_0)}{f_X(x_0)} = \frac{k-1}{x_*} - \frac{k x_*^{k-1}}{\l^k}.
\] 
Now, equation~(\ref{h_y_prime_x0_linear}) is equivalent to
\begin{equation}\label{case3} \begin{aligned}
\lim_{y\to\infty} &\Bigg(-\log^2 y \cdot \frac{\s_1\s_2\exp(-\s_2 x_2)}{2\s_0^2} \\
&\text{\hspace{1cm}}\ + \log y \cdot\left( \frac{(\mu_0 + \mu_1 x_2^{\mu_2}) \s_1\s_2 \exp\left(-\sigma_2 x_2\right)}{\s_0^2}+  \frac{\mu_1\mu_2 x_2^{\mu_2-1}}{\s_0} \right) \\
&\text{\hspace{1cm}}\ - \frac{(\mu_0 + \mu_1 x_2^{\mu_2})\mu_1\mu_2 x_2^{\mu_2-1}}{\s_0}  -\frac{\s_1\s_2 \exp(-\s_2 x_2)}{2\s_0}  \\
&\text{\hspace{1cm}}\ - \frac{ (\mu_0 + \mu_1 x_2^{\mu_2})^2 \cdot  \s_1\s_2 \exp\left(-\sigma_2 x_2\right) }{2\s_0^2}\\
&\text{\hspace{1cm}}\ + \frac{\mu_1 \mu_2 x_2^{\mu_2-1}}{\log y - \mu_0 -\mu_1x_2^{\mu_2}} + \frac{k-1}{x_2} - \frac{k x_2^{k-1}}{\l^k}\Bigg) = 0.
\end{aligned} \end{equation}
The dominating terms in equation~(\ref{case3}) are of the order $\log^2 y$, $\log y \cdot x_2^{\mu_2-1}$ and $x_2^{k-1}$. So, we can simplify equation~(\ref{case3}) to
\begin{equation}\label{case32} \begin{aligned}
\lim_{y\to\infty} -\log^2 y \cdot \frac{\s_1\s_2\exp(-\s_2 x_2)}{2\s_0^2} + \log y \cdot\frac{\mu_1\mu_2 x_2^{\mu_2-1}}{\s_0} - \frac{k x_2^{k-1}}{\l^k} = 0.
\end{aligned} \end{equation}
We note that the first and third terms have a negative sign, and the second has a positive sign. We note that we cannot simplify this further without considering the following two options as $y\to\infty$: (a) $\exp(-\s_2 x_2)\log^2 y \gg x_2^{k-1}$; (b) $\exp(-\s_2 x_2)\log^2 y \ll x_2^{k-1}$. Both of these cases will yield a solution to equation~(\ref{case32}) which we call $x_{2a}$ and $x_{2b}$ respectively.

\bigskip

\noi \textbf{Case (iii-a): $x_*=x_{2a}\to\infty$ and $\exp(-\s_2 x_{2a})\log^2 y \gg x_{2a}^{k-1}$}

\bigskip

\noi We derive from equation~(\ref{case32}) that $x_{2a}$ must satisfy as $y\to\infty$
\[
 -\log y \cdot \frac{\s_1\s_2}{2\s_0}\exp(-\s_2 x_{2a}) + \mu_1\mu_2 x_{2a}^{\mu_2-1} = O\left(\frac{x_{2a}^{k-1}}{\log y}\right).
\]
So, $x_{2a}$ must satisfy as $y\to\infty$
\[
x_{2a}^{\mu_2-1}\exp(\s_2 x_{2a}) = \log y \cdot \left(\frac{\s_1\s_2}{2\s_0 \mu_1\mu_2} +  O\left(\frac{x_{2a}^{k-1}}{\exp(-\s_2 x_{2a})\log^2 y}\right)\right). \]
Finally, we derive the following asymptotic expression for $x_{2a}$ as $y\to\infty$
\begin{equation}\label{2a}
x_{2a} = \frac{\log\log y }{\s_2}+ O(\log\log\log y).
\end{equation}

\bigskip

\noi \textbf{Case (iii-b): $x_*=x_{2b}\to\infty$ and $\exp(-\s_2 x_{2a})\log^2 y \ll x_{2a}^{k-1}$}

\bigskip

\noi We derive from equation~(\ref{case32}) that $x_{2b}$ must satisfy as $y\to\infty$
\[
\log y \cdot\frac{\mu_1\mu_2}{\s_0} - \frac{k x_2^{k-\mu_2}}{\l^k} = O\left(\log^2 y\exp(-\s_2 x_2) x_2^{1-\mu_2}\right).
\]
So, $x_{2b}$ must satisfy as $y\to\infty$
\begin{equation}\label{2b}
 x_{2b} = \left(\frac{ \l^k\mu_1\mu_2}{k \s_0}\right)^{\frac{1}{k-\mu_2}} \cdot  \left(\log y \right)^{\frac{1}{k-\mu_2}} + O\left((\log y)^{\frac{1}{k-\mu_2} - \frac{k-2\mu_2}{k-\mu_2}}\right).
\end{equation}

\subsection{Identifiying local maxima and local minima}
In the previous section, we have found expressions for local extrema, see equations~(\ref{0}),~(\ref{2a}) and~(\ref{2b}). In this section, we will show by using the second derivative $h_y''$ that $x_0$ and $x_{2b}$ correspond to local maxima and that $x_{2a}$ corresponds to a local minimum.

We calculated before
\[\begin{aligned}
h_y'(x)
& = \frac{-\phi\left(p_y(x)\right)}{\overline{\Phi}\left(p_y(x)\right)} \cdot p_y'(x) + \frac{f'_X(x)}{f_X(x)}.
\end{aligned}\]
So,
\[
\begin{aligned}
h_y''(x)
& = - \frac{\phi\left(p_y(x)\right)^2 p_y'(x)^2}{\overline{\Phi}\left(p_y(x)\right)^2}  - \frac{\phi'\left(p_y(x)\right)p_y'(x)^2}{\overline{\Phi}\left(p_y(x)\right)} - \frac{\phi\left(p_y(x)\right)p_y''(x)}{\overline{\Phi}\left(p_y(x)\right)} - \frac{f'_X(x)^2}{f_X(x)^2} + \frac{f''_X(x)}{f_X(x)}.
\end{aligned}\]
We can simplify $h_y''(x_*)$ for $x_*\in\{x_0,x_{2a},x_{2b}\}$ as $y\to\infty$ by using the identities $\phi(p_y(x_*))/\overline{\Phi}(p_y(x_*)) \sim p_y(x_*) + p_y(x_*)^{-1}$ as $y\to\infty$ and $\phi'(x) =-x\phi(x)$. We get as $y\to\infty$
\[\begin{aligned}
h_y''(x_*) &\sim  - \left(p_y(x_*)+\frac{1}{p_y(x_*)}\right)^2 p_y'(x_*)^2 +  \left(p_y(x_*)+\frac{1}{p_y(x_*)}\right) p_y(x_*) p_y'(x_*)^2 \\
&\text{\hspace{1cm}}\  - \left(p_y(x_*)+\frac{1}{p_y(x_*)}\right) p_y''(x_*) - \frac{f'_X(x_*)^2}{f_X(x_*)^2} + \frac{f''_X(x_*)}{f_X(x_*)} \\
&\sim  - p_y'(x_*)^2 - \frac{p_y'(x_*)^2 }{p_y(x_*)^2}  -  p_y(x_*)p_y''(x_*)-\frac{p_y''(x_*)}{p_y(x_*)}- \frac{f'_X(x_*)^2}{f_X(x_*)^2} + \frac{f''_X(x_*)}{f_X(x_*)} \\
&\sim  - p_y'(x_*)^2  -  p_y(x_*)p_y''(x_*)- \frac{f'_X(x_*)^2}{f_X(x_*)^2} + \frac{f''_X(x_*)}{f_X(x_*)}. \\
\end{aligned}\]
We work out $p_y'(x)^2$ and $p_y''(x)$ in terms of $p_y(x)$
\[
p_y'(x)^2 =\left(- p_y(x)\cdot\frac{\sigma'(x)}{\sigma(x)} - \frac{\mu'(x)}{\sigma(x)}\right)^2 = p_y(x)^2 \cdot  \frac{\sigma'(x)^2}{\sigma(x)^2} + p_y(x)\cdot \frac{2\sigma'(x)\mu'(x)}{\sigma(x)^2} + \frac{\mu'(x)^2}{\sigma(x)^2}
\]
and
\[
\begin{aligned} 
p_y''(x) &= \frac{\rd^2}{\rd x^2}\left(\frac{\log y - \mu(x)}{\sigma(x)}\right) \\
&= -\frac{\mu''(x)}{\sigma(x)} + 2\cdot \frac{\mu'(x)\sigma'(x)}{\sigma(x)^2} + (\log y - \mu(x)) \cdot\left( \frac{2\sigma'(x)^2}{\sigma(x)^3} - \frac{\sigma''(x)}{\sigma(x)^2}\right) \\
& = p_y(x) \cdot\left( \frac{2\sigma'(x)^2}{\sigma(x)^2} - \frac{\sigma''(x)}{\sigma(x)} \right)+ 2\cdot \frac{\mu'(x)\sigma'(x)}{\sigma(x)^2}-\frac{\mu''(x)}{\sigma(x)}.
\end{aligned}
\]
So, as $y\to\infty$
\[
\begin{aligned}
h_y''(x_*)
& \sim - p_y(x_*)^2 \cdot  \frac{\sigma'(x_*)^2}{\sigma(x_*)^2} - p_y(x_*)\cdot \frac{2\sigma'(x_*)\mu'(x_*)}{\sigma(x_*)^2} - \frac{\mu'(x_*)^2}{\sigma(x_*)^2}\\
&\text{\hspace{1cm}}\  -  p_y(x_*)\left(p_y(x_*) \cdot\left(\frac{2\sigma'(x_*)^2}{\sigma(x_*)^2} -  \frac{\sigma''(x_*)}{\sigma(x_*)} \right)+ 2\cdot \frac{\mu'(x_*)\sigma'(x_*)}{\sigma(x_*)^2}-\frac{\mu''(x_*)}{\sigma(x_*)}\right)\\
&\text{\hspace{1cm}}\ - \frac{f'_X(x_*)^2}{f_X(x_*)^2} + \frac{f''_X(x_*)}{f_X(x_*)} \\
& \sim - p_y(x_*)^2 \cdot\left(\frac{3\sigma'(x_*)^2}{\sigma(x_*)^2} - \frac{\sigma''(x_*)}{\sigma(x_*)}\right)- p_y(x_*)\left(\frac{4\mu'(x_*)\sigma'(x_*)}{\sigma(x_*)^2}-\frac{\mu''(x_*)}{\sigma(x_*)}\right) \\
&\text{\hspace{1cm}}\ - \frac{\mu'(x_*)^2}{\sigma(x_*)^2}- \frac{f'_X(x_*)^2}{f_X(x_*)^2} + \frac{f''_X(x_*)}{f_X(x_*)}. \\
\end{aligned}\]
For $x=x_0$, we have 
\[\begin{aligned}
\mu'(x_0) &= \mu_1\mu_2 x_0^{\mu_2-1}, \\
\mu''(x_0) &= -\mu_1 \mu_2(1-\mu_2)x_0^{\mu_2-2}, \\
\sigma(x_0) &\sim\sqrt{\s_0+\s_1}, \\
\sigma'(x_0) &\sim -\s_1\s_2/(2\sqrt{\s_0+\s_1}), \\
\sigma''(x_0)&\sim \sigma_1^2\sigma_2^2/(4(\s_0+\s_1)^{3/2}),\ \text{and} \\
p_y(x_0) &\sim \log y/\sqrt{\s_0+\s_1}.
\end{aligned}\]
So, 
\[\begin{aligned}
h_y''(x_0) 
&\sim- \frac{\log^2 y}{\s_0+\s_1} \cdot\left(\frac{3\s_1^2\s_2^2}{4(\s_0+\s_1)^2} - \frac{ \sigma_1^2\sigma_2^2}{4(\s_0+\s_1)^{2}}\right) \\
&\text{\hspace{1cm}}\ + \frac{\log y}{\sqrt{\s_0+\s_1}} \left(\frac{2 \mu_1\mu_2 x_0^{\mu_2-1} \cdot \s_1\s_2}{(\s_0+\s_1)^{3/2}}-\frac{\mu_1 \mu_2(1-\mu_2)x_0^{\mu_2-2}}{\sqrt{\s_0+\s_1}}\right) \\
&\text{\hspace{1cm}}\ - \frac{ \mu_1^2\mu_2^2 x_0^{2\mu_2-2}}{\s_0+\s_1}+ \frac{1}{x_0^2} + \frac{\log x_0 - \th}{x_0^2\a^2} - \frac{1}{x_0^2\a^2}  \\
&\sim-\frac{\s_1^2\s_2^2}{2(\s_0+\s_1)^3}\cdot \log^2 y  -  \frac{\mu_1 \mu_2(1-\mu_2)}{\s_0+\s_1} \cdot \log y \cdot x_0^{\mu_2-2}\\
&\text{\hspace{1cm}}\ - \frac{ \mu_1^2\mu_2^2}{\s_0+\s_1}\cdot  x_0^{2\mu_2-2} + \frac{1}{x_0^2} + \frac{\log x_0 - \th}{x_0^2\a^2} - \frac{1}{x_0^2\a^2}.  \\
\end{aligned}\]
We combine this result with equation~(\ref{0}), to get
\[\begin{aligned}
h_y''(x_0)
&\sim -\frac{\mu_1 \mu_2(1-\mu_2)}{\s_0+\s_1} \cdot \log y \cdot x_0^{\mu_2-2} \sim -C \left(\log y\right)^{2+\frac{1}{1-\mu_2}} \\
\end{aligned}\]
with
\[
C = \frac{\mu_1 \mu_2(1-\mu_2)}{\s_0+\s_1} \cdot  \left(\frac{\s_1 \s_2}{2\mu_1\mu_2(\s_0+\s_1)}\right)^{1+\frac{1}{1-\mu_2}}.
\]
We conclude that $h_y''(x_0)<0$ and that indeed $x_0$ corresponds to a local maximum.

For $x=x_{2*}$ with $*=a,b$, we have 
\[
\begin{aligned}
\mu'(x_{2*}) &= \mu_1\mu_2 x_{2*}^{\mu_2-1}, \\
\mu''(x_{2*}) &= -\mu_1\mu_2(1-\mu_2) x_{2*}^{\mu_2-2}, \\
\sigma(x_{2*}) &\sim \sqrt{\s_0}, \\
\sigma'(x_{2*}) &\sim-\s_1\s_2/(2\sqrt{\s_0}) \cdot \exp(-\s_2 x_{2*}), \\
\sigma''(x_{2*}) &\sim \s_1\s_2^2/(2\sqrt{\s_0}) \cdot \exp(-\s_2x_{2*}),\ \text{and} \\
p_y(x_{2b}) &\sim \log y/\sqrt{\s_0}.\end{aligned}\]
So,
\[\begin{aligned}
h_y''(x_{2*})
&\sim - \frac{\log^2 y}{\s_0} \cdot\left(\frac{3\s_1^2\s_2^2\cdot \exp(-2\s_2 x_{2*})}{4\s_0^2} - \frac{\s_1\s_2^2 \cdot \exp(-\s_2 x_{2*})}{2\s_0}\right) \\
&\text{\hspace{1cm}}\ +\frac{\log y}{\sqrt{\s_0}} \left(\frac{2\mu_1\mu_2 x_{2*}^{\mu_2-1} \cdot \s_1\s_2 \exp(-\s_2 x_{2*})}{\s_0^{3/2}}-\frac{\mu_1\mu_2(1-\mu_2) x_{2*}^{\mu_2-2}}{\sqrt{\s_0}}\right)  \\
&\text{\hspace{1cm}}\ - \frac{ \mu_1^2\mu_2^2 x_{2*}^{2\mu_2-2}}{\s_0}- \frac{k-1}{x_{2*}^2} - \frac{k(k-1) x_{2*}^{k-2}}{\l^k}.
\end{aligned}\]
For $x_{2*}=x_{2a}$, we simplify
\[
h_y''(x_{2a}) \sim \frac{\s_1\s_2^2}{2\s_0^2}\cdot \log^2 y \cdot \exp(-\s_2 x_{2a}) \\ 
\]
which confirms that $x_{2a}$ corresponds to a local minimum. Finally, for $x_{2*} = x_{2b}$, we simplify
\[\begin{aligned}
h_y''(x_{2b})&\sim - \frac{\mu_1\mu_2(1-\mu_2) x_{2b}^{\mu_2-2}}{\s_0}\cdot \log y - \frac{k(k-1)x_{2b}^{k-2}}{\l^k} \\
&\sim - \left(\frac{\mu_1\mu_2(1-\mu_2) \left(\frac{ \l^k\mu_1\mu_2}{k \s_0}\right)^{\frac{\mu_2-2}{k-\mu_2}} }{\s_0} + \frac{k(k-1) \left(\frac{ \l^k\mu_1\mu_2}{k \s_0}\right)^{\frac{k-2}{k-\mu_2}} }{\l^k}\right) \left(\log y \right)^{\frac{k-2}{k-\mu_2}} \\
&\sim - \left(\frac{\mu_1\mu_2(1-\mu_2) }{\s_0} \frac{k \s_0}{ \l^k\mu_1\mu_2} + \frac{k(k-1) }{\l^k}\right)\cdot \left(\frac{ \l^k\mu_1\mu_2}{k \s_0}\right)^{\frac{k-2}{k-\mu_2}}  \left(\log y \right)^{\frac{k-2}{k-\mu_2}}  \\
&\sim - \frac{k}{\l^k} \left(k-\mu_2\right)\cdot \left(\frac{ \l^k\mu_1\mu_2}{k \s_0}\right)^{\frac{k-2}{k-\mu_2}}  \left(\log y \right)^{\frac{k-2}{k-\mu_2}}.
\end{aligned}\]
Finally, it is clear to see that $h_y''(x_{2b})<0$ which confirms that $x_{2b}$ is a local maximum.

\subsection{Calculating the survival function of $Y$}

We will apply Proposition~\ref{laplace_approximation_variant4} to $g_y$ from equation~(\ref{gy22}), where we find that $k_0=2$ and $x_y^*=x_0(y)$. This gives us a lower bound for $\overline{F}_Y(y)$ as $y\to\infty$. We start with evaluating $g_y(x_0)$ and after, we check the smoothness assumption of the proposition for $k_0=2$. Finally, we derive an upper bound that is of the same order as the lower bound. Hence, we can combine the lower and upper bound to get an estimate for the rate of convergence to $0$ of $\overline{F}_Y(y)$.

Before, we evaluate $g_y(x_0)$ and $h_y''(x_0)$, we first simplify $p_y(x_0)$ and $p_y(x_0)$ as $y\to\infty$. We have
\[\begin{aligned}
p_y(x_0) &= \frac{\log y}{\sqrt{\s_0+\s_1}} - \frac{\mu_0}{\sqrt{\s_0+\s_1}} + O\left((\log y)^{-\frac{\mu_2}{1-\mu_2}}\right) \\
\end{aligned}\]
and
\[\begin{aligned}
\frac{1}{p_y(x_0)} & =\frac{\sqrt{\s_0+\s_1}}{\log y} + O\left((\log y)^{-2-\frac{\mu_2}{1-\mu_2}}\right)=\frac{\sqrt{\s_0+\s_1}}{\log y} + O\left((\log y)^{-2}\right).
\end{aligned}\]
So,
\[\begin{aligned}
g_y(x_0) &= \phi \left(\frac{\log y - \mu(x_0)}{\sigma(x_0)}\right) \cdot \left( \frac{\sigma(x_0)}{\log y - \mu(x_0)} + O\left(\left( \frac{\sigma(x_0)}{\log y - \mu(x_0)}\right)^{3}\right)\right)f_X(x_0)  \\
&= \frac{1}{\sqrt{2\pi}} \exp\left\{-\frac{1}{2} \left(\frac{\log y - \mu(x_0)}{\sigma(x_0)}\right)^2\right\}\\
&\text{\hspace{0.5cm}}\ \cdot \left( \frac{\sigma(x_0)}{\log y - \mu(x_0)} + O\left(\left( \frac{\sigma(x_0)}{\log y - \mu(x_0)}\right)^{3}\right)\right)\frac{ \exp\left\{-\frac{(\log x_0 - \th)^2}{2\a^2}\right\} }{\sqrt{2\pi} x_0 \alpha}\\
&= \frac{\sqrt{\s_0+\s_1}}{2\pi \alpha } \exp\left\{-\frac{1}{2} \left( \frac{\log y}{\sqrt{\s_0+\s_1}} - \frac{\mu_0}{\sqrt{\s_0+\s_1}} + O\left((\log y)^{-\frac{\mu_2}{1-\mu_2}}\right)\right)^2\right\} \\
&\text{\hspace{0.5cm}}\ \cdot \frac{\frac{1}{\log y} + O\left((\log y)^{-2}\right)}{\left(\frac{\s_1 \s_2}{2\mu_1\mu_2(\s_0+\s_1)}\right)^{-\frac{1}{1-\mu_2}} \cdot \left(\log y\right)^{-\frac{1}{1-\mu_2}} + O\left(\log^{-2}(y)\right)}\\
&\text{\hspace{0.5cm}}\ \cdot \exp\left\{-\frac{\left(\log \left[\left(\frac{\s_1 \s_2}{2\mu_1\mu_2(\s_0+\s_1)}\right)^{-\frac{1}{1-\mu_2}} \cdot \left(\log y\right)^{-\frac{1}{1-\mu_2}} + O\left(\log^{-2}(y)\right)\right] - \th\right)^2}{2\a^2}\right\}  \\
& = \exp\left\{-\frac{1}{2(\s_0+\s_1)} \left( \log^2 y - 2\mu_0 \log y + O\left((\log y)^{\frac{1-2\mu_2}{1-\mu_2}}\right)\right)\right\}. \\
\end{aligned}\]
Next, we check the assumptions of Proposition~\ref{laplace_approximation_variant4}. First, we clearly have $h_y'(x_0)(-h_y''(x_0))^{-1/2} = 0$. Secondly, we have one clearly dominating term in the second derivative of $h_y$ near $x_0$, so it is enough to show that
\[\begin{aligned}
\lim_{y\to \infty}&\frac{h_y''\left(x_0 + \frac{x}{\sqrt{-h_y''(x_0)}}\right) }{h_y''(x_0)} \\&= \lim_{x\to 0} \frac{p_y\left(x_0 + \frac{x}{\sqrt{-h_y''(x_0)}}\right)}{p_y(x_0)} \cdot \frac{\mu''\left(x_0 + \frac{x}{\sqrt{-h_y''(x_0)}}\right)}{\mu''\left(x_0\right)} \cdot \frac{\sigma(x_0)}{\sigma\left(x_0 + \frac{x}{\sqrt{-h_y''(x_0)}}\right)}
\end{aligned}\]
is equal to $1$ for any fixed $x$. Since $(-h_y''(x_0))^{-1/2} \ll x_0$ as $y\to\infty$ and since $p_y$ and $\s$ are differentiable at $0=\lim_{y\to\infty} x_0$, it is clear that the first and third term of the equation above tend to $1$. Since $\mu''(0)$ does not exist, we would need to work out the term involving the second derivative of $\mu$ more carefully. We get
\begin{equation}\label{muEqPrf}
 \frac{\mu''\left(x_0 + \frac{x}{\sqrt{-h_y''(x_0)}}\right)}{\mu''\left(x_0\right)} = \left(1 + \frac{x}{x_0\sqrt{-h_y''(x_0)}}\right)^{\mu_2-2}.
\end{equation}
We note that $x_0\sqrt{-h_y''(x_0)}$ is asymptotically equal to a constant times $(\log y)^{(1-2\mu_2)/(2-2\mu_2)}$. Since $\mu_2<0.5$, the second term within the brackets in equation~(\ref{muEqPrf}) tends to $0$ when $y\to\infty$. This yields that the right hand side of equation~(\ref{muEqPrf}) converges to $1$ as $y\to\infty$. This is enough to show the smoothness assumption of the proposition. We get that for any fixed $\tilde{x}>0$, there exists a constant $C_1(\tilde{x})$ such that
\begin{align}
\int_{0}^{\infty} g_y(x)\,\rd x
&\geq \int_{x_0 - \frac{\tilde{x}}{\sqrt{-h_y''(x_0)}}}^{x_0 + \frac{\tilde{x}}{\sqrt{-h_y''(x_0)}}} g_y(x)\,\rd x \nonumber \\
&\geq C_1(\tilde{x}) g_y(x_0)\cdot\frac{1}{\sqrt{-h_y''(x_0)}} \nonumber \\
(\text{as}\ y\to\infty)\ &= \exp\left\{-\frac{1}{2(\s_0+\s_1)} \left[ \log^2 y - 2\mu_0 \log y + O\left((\log y)^{\frac{1-2\mu_2}{1-\mu_2}}\right)\right]\right\}. \label{lowerbound_HW_distY}
\end{align}
Next, we evaluate $g_y(x_{2b})$ but
first we work out
\[
p_y(x_{2b}) = \frac{\log y}{\sqrt{\s_0}} + O\left((\log y)^{\frac{\mu_2}{k-\mu_2}}\right)
\]
and
\[\begin{aligned}
\frac{1}{p_y(x_0)} & =\frac{\sqrt{\s_0}}{\log y} + O\left((\log y)^{-2+\frac{\mu_2}{k-\mu_2}}\right).
\end{aligned}\]
So,
\[\begin{aligned}
g_y(x_{2b}) &=  \phi \left(\frac{\log y - \mu(x_{2b})}{\sigma(x_{2b})}\right) \left(\frac{\sigma(x_{2b})}{\log y - \mu(x_{2b})}+O\left(\frac{\sigma(x_{2b})^3}{(\log y - \mu(x_{2b}))^3}\right)\right) f_X(x_{2b}) \\
& = \phi\left( \frac{\log y}{\sqrt{\s_0}} + O\left((\log y)^{\frac{\mu_2}{k-\mu_2}}\right) \right)\cdot \frac{\sqrt{\s_0}}{\log y}\cdot\left(1+O\left((\log y)^{-1+\frac{\mu_2}{k-\mu_2}}\right)\right) \\
&\text{\hspace{1cm}}\ \cdot \frac{k}{\l^k} \left(\frac{ \l^k\mu_1\mu_2}{k \s_0}\right)^{\frac{k-1}{k-\mu_2}} (\log y)^{\frac{k-1}{k-\mu_2}} \left(1+ O\left((\log y)^{- \frac{k-2\mu_2}{k-\mu_2}}\right)\right)\\
&\text{\hspace{1cm}}\ \cdot\exp\left\{-\frac{\left(\frac{\s_1 \l^k\mu_1 \mu_2}{k \s_0^2 }\right)^{\frac{k}{k-\mu_2}} (\log y)^{\frac{k}{k-\mu_2}}}{\l^k}\left(1 + O((\log y)^{-\frac{k-2\mu_2}{k-\mu_2}}\right)\right\}\\
& = \exp\left\{-\frac{1}{2 \s_0} \log^2(y) + O\left( (\log y)^{\frac{k}{k-\mu_2}}\right)\right\}.
\end{aligned}\]
In particular, we find that $g_y(x_0) > g_y(x_{2b})$ for $y$ large enough. We have now all tools available to find an upperbound that gives the result directly,
\[
g_y(x) \leq \tilde{g}_y(x) := \begin{cases} \max\{g_y(x):\ x\in[0,x_2]\} & \text{for}\ 0\leq x\leq x_2, \\ f_X(x) & \text{for}\ x>x_2. \end{cases}
\]
Since $g_y(x_{2b}) \leq g_y(x_0)$ for $y$ large enough, we have derived that the maximum over the interval $[0,x_2]$ is attained at $x_0$. We here note that we do not need to show that $x_3$ and $x_4$ cannot exist as per definition, as they would clearly need to be larger than $x_2$ if they exist. So, as $y\to\infty$
\begin{align}
&\int_{0}^{\infty} g_y(x)\,\rd x \leq \int_0^{x_{2b}} g_y(x_0)\,\rd x + \int_{x_{2b}}^{\infty} f_X(x)\,\rd x = x_{2b} g_y(x_0) + \overline{F}_X(x_{2b}) \nonumber \\
& =\left(\frac{ \l^k\mu_1\mu_2}{k \s_0}\right)^{\frac{1}{k-\mu_2}} \left(1 + O\left((\log y)^{- \frac{k-2\mu_2}{k-\mu_2}}\right)\right) \left(\log y \right)^{\frac{1}{k-\mu_2}} \\
&\text{\hspace{1cm}}\ \cdot\exp\left\{-\frac{1}{2(\s_0+\s_1)} \left[ \log^2 y - 2\mu_0 \log y + O\left((\log y)^{\frac{1-2\mu_2}{1-\mu_2}}\right)\right]\right\}  \nonumber \\
&\text{\hspace{0.5cm}}\ +\exp\left\{- \l^{-k}\left(\frac{ \l^k\mu_1\mu_2}{k \s_0}\right)^{\frac{k}{k-\mu_2}} (\log y)^{\frac{k}{k-\mu_2}} \left[ 1 + O\left((\log y)^{-\frac{k-2\mu_2}{k-\mu_2}}\right)\right]\right\} \nonumber \\
& = \exp\left\{-\frac{1}{2(\s_0+\s_1)} \left[ \log^2 y - 2\mu_0 \log y + O\left((\log y)^{\frac{1-2\mu_2}{1-\mu_2}}\right)\right]\right\} .  \label{upperbound_HW_distY}
\end{align}
Now, combining equation~(\ref{lowerbound_HW_distY}) and equation~(\ref{upperbound_HW_distY}), yields as $y\to\infty$
\[
\P(Y>y) = \int_{0}^{\infty} g_y(x)\,\rd x =  \exp\left\{-\frac{1}{2(\s_0+\s_1)} \left[\log^2 y - 2\mu_0 \log y + O\left((\log y)^{\frac{1-2\mu_2}{1-\mu_2}}\right)\right]\right\}.
\]

\subsection{Calculating $\eta$}\label{hwmodeldetails_supmat_6}
We use the previous work to transform $Y$ to $Y_E$ on standard exponential margins. Thus
\[\begin{aligned}
Y_E &= F_E^{-1}(F_Y(Y)) = -\log(1-F_Y(Y)) \\
&= - \log\left( \exp\left\{-\frac{1}{2(\s_0+\s_1)} \left[ \log^2 Y - 2\mu_0 \log Y + O\left((\log Y)^{\frac{1-2\mu_2}{1-\mu_2}}\right)\right]\right\} \right) \\
&= \frac{1}{2(\s_0+\s_1)} \left( \log^2 Y - 2\mu_0 \log Y + O\left((\log Y)^{\frac{1-2\mu_2}{1-\mu_2}}\right)\right).
\end{aligned}\]
So, the function $T$ that transforms $\log Y$ to $Y_E$ is given by
\[
T(y) = \frac{y^2}{2(\s_0+\s_1)} - \frac{\mu_0 y}{\s_0+\s_1}+ O\left(y^{\frac{1-2\mu_2}{1-\mu_2}}\right),
\]
as $y\to\infty$. In calculating the extremal dependence measures, we need to solve $T(y)=u$ for large $y$. We get
\[
T^{-1}(u) = \sqrt{2(\s_0+\s_1)u} + O(1)
\]
as $u\to\infty$. We write down a formula for $\chi= \lim_{u\to\infty} \P(Y_E > u\mid (X/\l)^k > u)$ as $u\to\infty$
\begin{align}
&\P(Y_E > u\mid (X/\l)^k > u)  = e^{u}\int_{\l u^{1/k}}^{\infty} \P(T(\log Y) > u\mid X=x)f_X(x)\,\rd x \nonumber\\
& = e^{u}\int_{\l u^{1/k}}^{\infty} \P(\log Y > T^{-1}(u)\mid X=x)f_X(x)\,\rd x \nonumber\\
& = e^{u} \int_{\l u^{1/k}}^{\infty} \overline{\Phi}\left(\frac{T^{-1}(u) - \mu(x)}{\sigma(x)}\mid X=x\right)\cdot \frac{k x^{k-1}}{\l^k} \exp\left\{-\left(\frac{x}{\l}\right)^k\right\}\,\rd x. \nonumber
\end{align}
In particular, we have for $I=[\l u^{1/k},\l(2+\s_1/\s_0)^{1/k}]$
\begin{align}\label{chi_HW_integral}
\P(Y_E > u&\mid (X/\l)^k > u) \\
&> e^{u} \int_{I} \overline{\Phi}\left(\frac{T^{-1}(u) - \mu(x)}{\sigma(x)}\mid X=x\right)\cdot \frac{k x^{k-1}}{\l^k} \exp\left\{-\left(\frac{x}{\l}\right)^k\right\}\,\rd x.
\end{align}
For ease of presentation, we define $p_u(x) = [T^{-1}(u) - \mu(x)]/\sigma(x)$. Similar to the previous section, we define $g_u$ as the integrand and $h_u:=\log g_u$ as the log of the integrand, both are specified only on the integration domain $I$. For $x$ in the integration domain, we have
\[\begin{aligned}
h_u(x) &:= \log\left(\overline{\Phi}(p_u(x))f_X(x)\right).\\
\end{aligned}\]
We apply Proposition~\ref{laplace_approximation_variant4} to bound integral~(\ref{chi_HW_integral}) from below. In particular, we first need to find the mode of $h_u$ over the integration domain. Let $x_u$ be a sequence such that for each $u$, $x_u$ lies in the integration domain. So, then we can write $x = C_u u^{1/k} + o(u^{1/k})$ for some bounded set of constants $C_u \in [\l,  \l(2+\s_1/\s_0)]$. We have
\[\begin{aligned}
h_u'(x) &= -\frac{\phi(p_u(x))}{\overline{\Phi}(p_u(x))} \cdot p_u'(x) -  \frac{k-1}{x} - \frac{k x^{k-1}}{\l^{k}} \\
& = \frac{\phi(p_u(x))}{\overline{\Phi}(p_u(x))} \cdot \left(p_u(x)\cdot\frac{\sigma'(x)}{\sigma(x)} + \frac{\mu'(x)}{\sigma(x)}\right) -  \frac{k-1}{x} - \frac{k x^{k-1}}{\l^{k}}.
\end{aligned}\]
Since, $p_u(x) \sim \sqrt{2(1+\s_1/\s_0) u} \to\infty$ as $u\to\infty$, we simplify
\[\begin{aligned}
h_u'(x) 
&\sim \sqrt{2\left(1+\frac{\s_1}{\s_0}\right) u}\cdot \Bigg(\sqrt{2\left(1+\frac{\s_1}{\s_0}\right) u}\cdot\frac{-\s_1\s_2 e^{-\s_2\left( C_u u^{1/k} + o(u^{1/k})\right)}}{2\s_0} \\
&\text{\hspace{5cm}}\ + \frac{\mu_1\mu_2 \left( C_u u^{1/k} + o(u^{1/k})\right)^{\mu_2-1}}{\sqrt{\s_0}}\Bigg) \\
&\text{\hspace{1cm}}\ -  \frac{k-1}{ C_u u^{1/k} + o(u^{1/k})} - \frac{k \left( C_u u^{1/k} + o(u^{1/k})\right)^{k-1}}{\l^{k}}.  \\
&\sim -\frac{k C_u^{k-1} u^{1-1/k}}{\l^{k}}.
\end{aligned}\]
In particular, we derive that $h_u'(x)<0$ as $u\to\infty$. So, the maximum of $h_u$ over the integration domain must be attained at the boundary and hence is given by $x_0=\l u^{1/k}$. In particular, we get $h_u'(x_0) \sim -k u^{1-1/k}/\l$. We now will show that we can apply Proposition~\ref{laplace_approximation_variant4} with $k_0=1$. We have, as $u\to\infty$, 
\[\begin{aligned}
h_u(\l u^{1/k})
& =  -\frac{1}{2}\log(2\pi) - \frac{1}{2}p_u(\l u^{1/k})^2 - \log p_u(\l u^{1/k}) + \log f_X(\l u^{1/k})\\
& =  -\frac{1}{2}\log(2\pi) - \frac{1}{2}\left(\frac{T^{-1}(u)  - \mu(\l u^{1/k})}{\s(\l u^{1/k})}\right)^2  \\
&\text{\hspace{1cm}}\ - \log \left(\frac{T^{-1}(u) - \mu(\l u^{1/k})}{\s(\l u^{1/k})}\right) + \log \left(\frac{k u^{(k-1)/k} }{\l}\right) -u\\
& = - \left(2+\frac{\s_1}{\s_0}\right) u + O\left(u^{1/2 + \mu_2/k}\right). \\
\end{aligned}\]

Next, we check the smoothness assumption of Proposition~\ref{laplace_approximation_variant4} with $k_0=1$. Let $\d>0$ and $0\leq x\leq \d$. It is now enough to show that the limit of $u$ to infinity of the following expression tends to $1$. We have
\[
\lim_{u\to\infty}\frac{h_u'\left(\l u^{1/k} + \frac{x}{-h_u'(\l u^{1/k})}\right) }{h_u'(\l u^{1/k})} 
= \lim_{u\to\infty} \frac{\left(\l u^{1/k} + \frac{\l x}{k u^{1-1/k}}\right)^{k-1}}{ u^{(k-1)/k}}
= \lim_{u\to\infty} \left(\l  + \frac{\l x}{k }\right)^{k-1} = 1.
\]
This is enough to show the smoothness assumption of Proposition~\ref{laplace_approximation_variant4} with $k_0=1$. We conclude that for each $\tilde{x}$, there exists a constant $C_1(\tilde{x})$ such that
\begin{align}
\int_{\l u^{1/k}}^{\infty} g_u(x)\,\rd x 
&\geq \int_{\l u^{1/k}}^{\l u^{1/k}+\frac{\tilde{x}}{-h_u'(x_0)}} g_u(x)\,\rd x \geq C_1(\tilde{x}) g_u(\l u^{1/k}) \cdot \frac{1}{-h_u'(\l u^{1/k})} \nonumber \\
&\stackrel{(\text{as}\ u\to\infty)}{=} e^{-\left(2  +\frac{\s_1}{\s_0}\right) u + O\left(u^{\frac{1}{2} + \frac{\mu_2}{k}}\right)}. \label{lowerbound_chi_HW}
\end{align}
To get an upper bound, we use the following crude upper bound $\tilde{g}_u$ for $g_u$,
\[
g_u(x) \leq \tilde{g}_u(x):=\begin{cases} g_u(\l u^{1/k}) &\text{for}\ \l u^{1/k} \leq x \leq \l \left(2  +\frac{\s_1}{\s_0}\right)^{1/k} u^{1/k}, \\ f_X(x) &\text{for}\ x > \l \left(2  +\frac{\s_1}{\s_0}\right)^{1/k}  u^{1/k}.\end{cases}
\]
We get as $u\to\infty$,
\begin{align}
\int_{\l u^{1/k}}^{\infty} g_u(x)\,\rd x 
&\leq \left(\l \left(2  +\frac{\s_1}{\s_0}\right)^{1/k} u^{1/k} - \l u^{1/k}\right) g_u(\l u^{1/k}) \\
&\text{\hspace{1cm}}\ + \overline{F}_X\left(\l \left(2  +\frac{\s_1}{\s_0}\right)^{1/k}u^{1/k}\right) \nonumber \\
&= \exp\left(- \left(2  +\frac{\s_1}{\s_0}\right) u + O\left(u^{1/2 + \mu_2/k}\right) \right) + \exp\left(-\left(2  +\frac{\s_1}{\s_0}\right)u\right) \nonumber\\
&= \exp\left(- \left(2  +\frac{\s_1}{\s_0}\right) u + O\left(u^{1/2 + \mu_2/k}\right) \right). \label{upperbound_chi_HW}
\end{align}
Combining equations~(\ref{lowerbound_chi_HW}) and~(\ref{upperbound_chi_HW}), we get
\[
\P(Y_E > u\mid (X/\l)^k > u) = \int_{\l u^{1/k}}^{\infty} g_u(x)\,\rd x = \exp\left(-  \left(2  +\frac{\s_1}{\s_0}\right) u + O\left(u^{1/2 + \mu_2/k}\right) \right)
\]
as $u\to\infty$.
From this expression, it is straightforward to see that $\xi=0$ and
\[
\eta^{-1} =2 +\frac{\s_1}{\s_0}.
\]

\section{Details on Calculations for the Exact HT model}\label{supHT_section}
\subsection{Introduction}
Assume model~(\ref{HT_formulation}) for random vector $(X,Y)$ with $\overline{H}$ as in equation~(\ref{H_parameterisation}). We recall that $(X,Y)$ is a random vector such that $X$ and $Y$ both have standard Laplace margins. Moreover,there exist $0\leq \a\leq 1$, $\b<1$ and $u>0$ such that for $x>u$
\[
\P(Y> y\mid X=x) = \overline{H}\left(\frac{y - \a x}{x^{\b}}\right),
\]
holds for all $y\in\R$ with
\[
\overline{H}\left(z\right) = \exp(- \gamma z^{\delta})\1\{z>0\} + \1\{z\leq 0\}
\]
for $\gamma>0$ and $\delta\geq(1-\b)^{-1}$. In this section, we work out the value for $\eta$ when $0<\alpha<1,\ \beta>0$ and $\d=(1-\b)^{-1}$. The other cases are significantly easier to work out and the results of these cases are stated in the main paper.

\subsection{Calculating $\eta$} 
We write
\[\begin{aligned}
\P(Y>u,X>u) &= \int_u^{\infty}  \overline{H}\left(\frac{u - \a x}{x^{\b}}\right) f_X(x)\,\rd x \\
&= \frac{1}{2}\int_u^{u/\alpha}  \exp\left(-\g\left(\frac{u - \a x}{x^{\b}}\right)^{\delta} - x\right)\,\rd x + \frac{1}{2}\int_{u/\alpha}^{\infty} \exp(-x)\,\rd x\\
&= \frac{1}{2}\int_u^{u/\alpha}  \exp\left(-\g\left(\frac{u - \a x}{x^{\b}}\right)^{\delta} - x\right)\,\rd x + \frac{1}{2}\exp\left(-\frac{u}{\alpha}\right).
\end{aligned}\]
In general, we cannot evaluate the first integral in closed form for finite $u$. However, we can bound it from below using Proposition~\ref{laplace_approximation_variant4}. A bound from above can again be found directly. We define the integration domain $I=[u,u/\alpha]$, 
\[
g_u(x) := \exp\left(-\g\left(\frac{u - \a x}{x^{\b}}\right)^{\delta} - x\right)
\]
for $x\in I$ and $h_u:=\log g_u$ on $I$. We now need to determine whether or not the mode $x_0:=x_0(u)$ of the integrand $g_u$ over the integration domain $I$ lies on the boundary of $I$ or in the interior of $I$. We assume that $x_0$ lies in the interior of $I$, then we have
\[\begin{aligned}
0=h_u'(x_0) & =\gamma \delta \left(\frac{u-\alpha x_0}{x_0^{\b}}\right)^{\delta-1} \cdot \left(\frac{\a}{x_0^{\b}} + \frac{(u-\alpha x_0)\b}{x_0^{\b+1}}\right) - 1 \\
&= \gamma \delta \beta (u-\a x_0)^{\d} x_0^{-\b\d-1}  +  \gamma \delta \a (u-\a x_0)^{\d-1} x_0^{-\b\d} - 1\\
&= \gamma \delta \beta (u-\a x_0)^{\d} x_0^{-\d}  +  \gamma \delta \a (u-\a x_0)^{\d-1} x_0^{-\d+1} - 1
\end{aligned}\]
and we derive that
\begin{equation}\label{eqn}
\beta (u-\a x_0)^{\d}  + \a (u-\a x_0)^{\d-1} x_0 = \frac{1}{\gamma\delta} x_0^{\d}.
\end{equation}
Since, we work under the premise that $x_0\in(u,u/\alpha)$, we are only interested in finding solutions that satisfy $x_0=\tilde{c} u + o(u)$ as $u\to\infty$ for some $\tilde{c}\in[1,1/\alpha]$, otherwise the mode of $h_u$ is found at the boundary of the integration domain at $u$. We try $x_0 = c u$ with $c\in(0,\infty)$ in equation~(\ref{eqn}), and we derive that this is an exact solution if $c$ solves
\begin{equation}\label{eq_for_c}
0 = \gamma\delta\left(\beta(1-\alpha c)^{\delta} + \alpha c(1-\alpha c)^{\delta-1} \right) - c^{\delta} = \gamma(1-\alpha c)^{\delta-1}\left(\d-1 +\alpha c\right) - c^{\delta}.
\end{equation}
Since the right hand side is a continuous function of $c$ for $c\in[0,1/\alpha]$, we show by the intermediate value theorem that $c\in(0,1/\alpha)$ by inserting $c=0$ and $c=1/\alpha$ and comparing signs of the right hand side of equation~(\ref{eq_for_c}). Indeed, for $c=0$, we have that
\[
\gamma(1-\alpha c)^{\delta-1}\left(\d-1 +\alpha c\right) - c^{\delta} = \gamma(\delta-1) > 0
\]
and for $c=1/\alpha$, we have that
\[
\gamma(1-\alpha c)^{\delta-1}\left(\d-1 +\alpha c\right) - c^{\delta} = -\alpha^{-\d}<0.
\]
We recall that we are only interested in the value for $c$ if $c\in(1,1/\alpha)$. Hence, let $c=1$ in the right hand side of equation~(\ref{eq_for_c}) to give
\begin{align}
\gamma\delta\left(\beta(1-\alpha c)^{\delta} + \alpha c (1-\alpha c)^{\delta-1} \right) - c^{\delta}
&=  \gamma(1-\alpha)^{\d-1}\left(\left(\d-1\right)(1-\alpha ) + \d\alpha  \right) - 1 \nonumber\\
& = \g(1-\a)^{\d-1}(\d-1+\a) -1\nonumber
\end{align}
which is negative if and only if $ \g(1-\a)^{\d-1}(\d-1+\a)<1$. We conclude that $c\in(0,1)$ if and only if $\g(1-\a)^{\d-1}(\d-1+\a)<1$ and $c\in[1,1/\alpha)$ if and only if $\g(1-\a)^{\d-1}(\d-1+\a)\geq1$. We term these cases as Case (2a) and Case (2b), respectively. In Case (2b), $x_0$ lies in the interior of the integration domain $I$ for large enough $u$, and in Case (2a), the mode over the integration domain $I$ is found at $u$ on the boundary.

We work out $g_u(x_0)$ for both Case (2a) and (2b),
\[\begin{aligned}
g_u(x_0)
& =   \exp\left\{- \gamma \left(\frac{u - \a x_0}{x_0^{\b}}\right)^{\delta} - x_0\right\} \\
& =   \exp\left\{- \gamma \left(\frac{u - \a (cu + o(u))}{(cu + o(u))^{\b}}\right)^{\delta} - c u  + o(u)\right\} \\
& =   \exp\left\{- \gamma u^{\delta-\b\d} \frac{(1 - \a c + o(1))^{\delta}}{c^{\b\d} + o(1))} - c u  + o(u)\right\} \\
& =  \exp\left\{- \left( \frac{\gamma (1 - \a c)^{\d}}{c^{\b\d}} + c \right)u  + o(u)\right\}. \\
\end{aligned} \]
Next, we work out $h_u'(x_0)$ in Case (2a)
\[\begin{aligned}
h_u'(x_0) 
& =  \gamma \delta \beta (u-\a u)^{\d} u^{-\b\d-1}  +  \gamma \delta \a (u-\a u)^{\d-1} u^{-\b\d} - 1 \\
& =  \gamma (1-\a)^{\d-1} \left(\delta-1+\a\right) - 1. \\
\end{aligned} \]
By definition of Case (2a), we have that $h_u'(x_0)<0$. Let $C>0$ and $|x|\leq C$, then as $u\to\infty$
\[\begin{aligned}
h_u'\left(x_0 - \frac{x}{h_u'(x_0)}\right)
& =\gamma \delta \beta \left(u-\a \left(x_0 + \frac{x}{-h_u'(x_0)}\right)\right)^{\d} \left(x_0 + \frac{x}{-h_u'(x_0)}\right)^{-\b\d-1}  \\
&\text{\hspace{0.5cm}}\ +  \gamma \delta \a \left(u-\a \left(x_0 + \frac{x}{-h_u'(x_0)}\right)\right)^{\d-1} \left(x_0 + \frac{x}{-h_u'(x_0)}\right)^{-\b\d}- 1 \\
& = \gamma \delta \beta \left(u-\a u- \frac{x\alpha}{1 - \gamma (1-\a)^{\d-1} \left(\delta-1+\a\right) }\right)^{\d}  \\
&\text{\hspace{1cm}}\ \cdot \left(u+ \frac{x}{1 - \gamma (1-\a)^{\d-1} \left(\delta-1+\a\right) }\right)^{-\b\d-1}  \\
&\text{\hspace{0.5cm}}\ +  \gamma \delta \a \left(u-\a u - \frac{x \alpha}{1 - \gamma (1-\a)^{\d-1} \left(\delta-1+\a\right) }\right)^{\d-1}  \\
&\text{\hspace{1cm}}\ \cdot \left(u + \frac{x}{1 - \gamma (1-\a)^{\d-1} \left(\delta-1+\a\right) }\right)^{-\b\d} - 1  \\
& = \gamma \delta \beta \left(u^{\delta}(1-\a)^{\delta}  + O\left(u^{\delta-1}\right)\right)  \left(u^{-\b\d-1} + O\left(u^{-\b\d-2}\right)\right)  \\
&\text{\hspace{0.5cm}}\ +  \gamma \delta \a \left(u^{\delta-1} (1-\a)^{\delta-1}+ O\left(u^{\delta-2}\right)\right)\left(u^{-\b\d} + O\left(u^{-\b\d-1}\right)\right)- 1  \\
& = h_u'(x_0) + O\left(u^{\delta-2-\b\d}\right).  \\
\end{aligned}\]
So,
\[
\lim_{u\to\infty}\frac{h_u'\left(x_0 + \frac{x}{-h_u'(x_0)}\right) }{h_u'(x_0)} = 1,
\]
which is enough to show the smoothness assumption of Proposition~\ref{laplace_approximation_variant4} with $k_0=1$. We get that for any fixed $\tilde{x}>0$ there exist a $C_1(\tilde{x})>0$ such that
\[
\begin{aligned}
\P(X>u,\ Y>u)
& = \int_{u}^{u/\alpha} g_u(x)\,\rd x + \frac{1}{2}\exp\left\{-\frac{u}{\alpha}\right\} \\
& \geq \frac{1}{2}\int_{x_0 - \frac{\tilde{x}}{-h_u'(x_0)}}^{x_0 + \frac{\tilde{x}}{-h_u'(x_0)}} g_u(x)\,\rd x+ \frac{1}{2}\exp\left\{-\frac{u}{\alpha}\right\} \\
& \geq \frac{1}{2} C_1(\tilde{x}) g_u(x_0) \cdot \frac{1}{-h_u'(x_0)}+ \frac{1}{2}\exp\left\{-\frac{u}{\alpha}\right\} \\
& \geq \frac{1}{2} C_1(\tilde{x}) \exp\left\{- \left( \gamma (1 - \a)^{\d} +1 \right)u  + o(u)\right\} \\
&\text{\hspace{1cm}}\ \cdot \frac{1}{1- \gamma (1-\a)^{\d-1} \left(\delta-1+\a\right) }+ \frac{1}{2}\exp\left\{-\frac{u}{\alpha}\right\} \\
& = \exp\left\{- \left(\gamma (1 - \a)^{\d} +1  \right)u  + o(u)\right\}. \\
\end{aligned}
\]
In the last step we used that $\gamma (1 - \a)^{\d} +1  < 1/\alpha$ holds, which can be directly derived from the assumptions corresponding to Case (2b). Similarly to before, we can find an upper bound rather straightforwardly using the following upperbound for $g_u(x)$
\[
g_u(x) \leq \tilde{g}_u(x):=\begin{cases} g_u(x_0) &\text{for}\ u \leq x \leq u/\alpha, \\ f_X(x) &\text{for}\ x > u/\alpha.\end{cases}
\]
So,
\[
\begin{aligned}
\P(X>u,\ Y>u)
& = \int_{u}^{u/\alpha} g_u(x)\,\rd x + \frac{1}{2}\exp\left\{-\frac{u}{\alpha}\right\} \\
& \leq u\left(\frac{1}{\alpha} - 1\right) g_u(x_0) + \frac{1}{2}\exp\left\{-\frac{u}{\alpha}\right\} \\
& = u\left(\frac{1}{\alpha} - 1\right)\exp\left\{- \left( \gamma (1 - \a)^{\d} +1  \right)u  + o(u)\right\} + \frac{1}{2}\exp\left\{-\frac{u}{\alpha}\right\} \\
& = \exp\left\{- \left(\gamma (1 - \a)^{\d} +1  \right)u  + o(u)\right\}. \\
\end{aligned}
\]
We conclude that
\[
\P(X>u,\ Y>u) = \exp\left\{- \left( \gamma (1 - \a)^{\d} +1  \right)u  + o(u)\right\},
\]
$\chi=0$ and
\[
\eta = \left(\gamma (1 - \a)^{\d} +1  \right)^{-1}.
\]
For Case (2a), we work out $h_u''(x_0)$ as $u\to\infty$
\[\begin{aligned}
h_u''(x_0)
& =  - \a^2\gamma \delta(\d-1)(u-\a x_0)^{\d-2} x_0^{-\b\d} -  2 \a \b \gamma \delta^2(u-\a x_0)^{\d-1} x_0^{-\b\d-1}   \\
&\text{\hspace{1cm}}\ -  \beta\delta(\b\d+1)\gamma  (u-\a x_0)^{\d} x_0^{-\b\d-2}  \\
& =  - \a^2\gamma \delta(\d-1)(u-\a (cu + o(u)))^{\d-2} (cu + o(u))^{-\b\d} \\
&\text{\hspace{1cm}}\ -  2 \a \b \gamma \delta^2(u-\a (cu + o(u)))^{\d-1} (cu + o(u))^{-\b\d-1}  \\
&\text{\hspace{1cm}}\ -  \beta\delta(\b\d+1)\gamma  (u-\a (cu + o(u)))^{\d} (cu + o(u))^{-\b\d-2}  \\
& =  - \Big[\a^2\gamma \delta(\d-1) (1-\a c)^{\d-2} c^{-\b\d} + 2 \a \b \gamma \delta^2  (1-\a c)^{\d-1} c^{-\b\d-1}   \\
&\text{\hspace{1cm}}\ + \beta\delta(\b\d+1)\gamma (1-\a c)^{\d} c^{-\b\d-2} \Big]u^{\d-2-\b\d} + o\left(u^{\d-2-\b\d}\right)\\
& =  - \beta\delta^2 \gamma c^{-\b\d-2} (1-\a c)^{\d-2} \left(\a^2  c^{2} + 2 \a   (1-\a c)c + (1-\a c)^{2} \right)u^{-1} + o\left(u^{-1}\right)\\
& =  - \beta\delta^2 \gamma c^{-\b\d-2} (1-\a c)^{\d-2} u^{-1} + o\left(u^{-1}\right)\\
& =  - \delta(\delta-1) \gamma c^{-\delta-1} (1-\a c)^{\d-2} u^{-1} + o\left(u^{-1}\right).\\
\end{aligned} \]
Now, let $C>0$ and $|x|\leq C$, then we have $x_0 + x(-h_u''(x_0))^{-1/2}= c u + o(u)$. So,
\[\begin{aligned}
h_u''\left(x_0+\frac{x}{\sqrt{-h_u''(x_0)}}\right)
& = h_u''\left(cu + o(u) \right). \\
\end{aligned} \]
So,
\[\lim_{u\to\infty} \frac{h_u''\left(x_0+\frac{x}{\sqrt{-h_u''(x_0)}}\right)}{h_u''(x_0)} =\lim_{u\to\infty}\frac{- \delta(\delta-1) \gamma c^{-\delta-1} (1-\a c)^{\d-2} u^{-1}(1+o(1))}{- \delta(\delta-1) \gamma c^{-\delta-1} (1-\a c)^{\d-2} u^{-1}(1+o(1))} = 1,\]
which is enough to show the smoothness assumption of Proposition~\ref{laplace_approximation_variant4} with $k_0=1$. We get that for any fixed $\tilde{x}>0$ there exist a $C_1(\tilde{x})>0$ such that as $u\to\infty$
\[\begin{aligned}
\P(X>u,Y>u) 
& = \int_{u}^{u/\alpha} g_u(x)\,\rd x + \frac{1}{2}\exp\left\{-\frac{u}{\alpha}\right\} \\
& \geq \frac{1}{2}\int_{x_0 - \frac{\tilde{x}}{-h_u'(x_0)}}^{x_0 + \frac{\tilde{x}}{-h_u'(x_0)}} g_u(x)\,\rd x  \\
& \geq \frac{1}{2} C_1(\tilde{x}) g_u(x_0) \cdot \frac{1}{-h_u'(x_0)} \\
& =  \exp\left\{-\left(\frac{\g(1-\a c)^{\d}}{c^{\d-1}} + c\right)u + o(u)\right\}.
\end{aligned}\]
Similarly to before, we can find an upper bound rather straightforwardly,
\[
\begin{aligned}
\P(X>u,\ Y>u)
& = \int_{u}^{u/\alpha} g_u(x)\,\rd x + \frac{1}{2}\exp\left\{-\frac{u}{\alpha}\right\} \\
& \leq u\left(\frac{1}{\alpha} - 1\right) g_u(x_0) + \frac{1}{2}\exp\left\{-\frac{u}{\alpha}\right\} \\
& = u\left(\frac{1}{\alpha} - 1\right) \exp\left\{-\left(\frac{\g(1-\a c)^{\d}}{c^{\b\d}} + c\right)u + o(u)\right\} + \frac{1}{2}\exp\left\{-\frac{u}{\alpha}\right\} \\
& =\exp\left\{-\left(\frac{\g(1-\a c)^{\d}}{c^{\d-1}} + c\right)u + o(u)\right\}.
\end{aligned}
\]
So,
\[
\P(X>u,\ Y>u) =\exp\left\{-\left(\frac{\g(1-\a c)^{\d}}{c^{\b\d}} + c\right)u + o(u)\right\},
\]
and we conclude that $\chi=0$ and 
\[
\eta = \left(\frac{\g(1-\a c)^{\d}}{c^{\d-1}} + c\right)^{-1}.
\]

\bibliography{bibpaper3}
\bibliographystyle{apa}
